\documentclass[11pt,a4paper]{article}
\usepackage{titlesec}
\usepackage{fancyhdr}
\usepackage{a4wide}
\usepackage{graphicx}
\usepackage{float}
\usepackage{amssymb}
\usepackage{amsmath}
\usepackage{amsthm}
\usepackage{color}
\usepackage{mathrsfs}
\usepackage{array}
\usepackage{multido}
\usepackage{wrapfig}
\usepackage[T1]{fontenc}
\usepackage{inputenc}
\usepackage[english]{babel}

\usepackage{lmodern}
\usepackage{hyperref}
\usepackage{geometry}
\usepackage[numbers]{natbib}
\hypersetup{
pdfpagemode=UseThumbs,
pdftoolbar=true,        % show Acrobat's toolbar?
pdfmenubar=true,        % show Acrobats menu?
pdffitwindow=false,     % window fit to page when opened
pdfstartview={Fit},    % fits the width of the page to the window
pdftitle={A numerical approach for the Poisson equation in a planar domain with a small inclusion},    % title
pdfauthor={L. Chesnel, X. Claeys},     % author
pdfsubject={},  % subject of the document
pdfcreator={L. Chesnel, X. Claeys},   % creator of the document
pdfproducer={L. Chesnel, X. Claeys}, % producer of the document
pdfkeywords={}, % list of keywords
pdfnewwindow=true,      % links in new window
colorlinks=true,       % false: boxed links; true: colored links
linkcolor=magenta,          % color of internal links
citecolor=red,        % color of links to bibliography
filecolor=cyan,      % color of file links
urlcolor=blue           % color of external links
}

\newcommand{\dsp}{\displaystyle}

\newcommand{\eps}{\varepsilon}
\newcommand{\om}{\omega}
\newcommand{\Om}{\Omega}
\newcommand{\mrm}[1]{\mathrm{#1}}

\newcommand{\bfx}{\boldsymbol{x}}
\newcommand{\bfs}{\mathbf{s}}
\newcommand{\bfxi}{\boldsymbol{\xi}}

\newcommand{\R}{\mathbb{R}}

\newcommand{\mL}{\mrm{L}}
\newcommand{\mH}{\mrm{H}}
\newcommand{\mV}{\mrm{V}}

\newcommand{\mD}{\mrm{D}}

\newcommand{\supp}{\mrm{supp}}

\newtheorem{lemma}{Lemma}[section]
\newtheorem{remark}{Remark}[section]

\newtheorem{proposition}{Proposition}[section]

\begin{document}

~\vspace{0.5cm}
\begin{center}
{\sc \bf\LARGE 
A numerical approach for the Poisson equation \\[6pt]  in a planar domain with a small inclusion}
\end{center}

\begin{center}
\textsc{Lucas Chesnel}$^1$, \textsc{Xavier Claeys}$^2$\\[16pt]
\begin{minipage}{0.9\textwidth}
{\small
$^1$ Centre de mathématiques appliquées, \'{E}cole Polytechnique, 91128 Palaiseau, France; \\
$^2$ Laboratoire Jacques-Louis Lions, Universit\'e Pierre et Marie Curie, 4 place Jussieu, 75005 Paris, France.\\[10pt] 
E-mail: \texttt{lucas.chesnel@cmap.polytechnique.fr}, \texttt{claeys@ann.jussieu.fr}\\[-14pt]
\begin{center}
(\today)
\end{center}
}
\end{minipage}
\end{center}
\vspace{0.4cm}

\noindent\textbf{Abstract.} 
We consider the Poisson equation in a domain with a small hole of size $\delta$. We present a simple numerical method, 
based on an asymptotic analysis, which allows to approximate robustly the far field of the solution as $\delta$ 
goes to zero without meshing the small hole. We prove the stability of the scheme and provide error estimates. We 
end the paper with numerical experiments illustrating the efficiency of the technique.\\

\noindent\textbf{Key words.} Small hole, asymptotic analysis, singular perturbation, finite element method.

\section{Introduction}
In the present article, we consider $\omega,\Omega\subset\R^{2}$ two bounded  Lipschitz domains 
such that $\overline{\omega}\subset\Omega$. Define $\omega_{\delta}:=\{\bfx\in \R^{2},\;\bfx/\delta\in \omega\}$ 
and $\Omega_{\delta}:=\Omega\setminus\overline{\omega}_{\delta}$. Given a data $f\in\mL^{2}(\Omega)$, 
we are interested in devising a robust and accurate numerical method for approximating, for small values of $\delta$, 
the \textit{far field} of the function satisfying  
\begin{equation}\label{PbInit}
u_{\delta}\in \mH^{1}_{0}(\Omega_{\delta})
\quad\textrm{ and }\quad
-\Delta u_{\delta} = f\quad\textrm{in}\;\Omega_{\delta}.
\end{equation}
In (\ref{PbInit}), $\mrm{H}^1_0(\Om_{\delta})$ denotes the subspace of the elements of the Sobolev 
space $\mrm{H}^1(\Om_{\delta})$ vanishing on $\partial\Om_{\delta}$. On the other hand, we call far field 
of $u_{\delta}$ the restriction of $u_{\delta}$ to $\Om\setminus\overline{\mD}_{r}$, where $\mD_{r} := \mD(0,r)$ is the disk with fixed arbitrary radius $r>0$. Problem (\ref{PbInit}), or variants of 
it, arises as a simple but relevant model in many applications ranging from electrical engineering 
\cite{865332,1353491} to flow transport around wells \cite{MR1990198,Peaceman}. This kind of problem 
also appears when considering wave scattering by small impenetrable inclusions \cite{MR2991243}.\\ 
\newline
In order to solve numerically Problem (\ref{PbInit}), a crude but rather natural idea would consist 
in neglecting the influence of the small inclusion on the total field $u_{\delta}$. Indeed (see for 
example \cite{MaNP00}), as $\delta\to 0$ the function $u_{\delta}$ converges toward
$u_{0}$ the solution to the limit problem where the inclusion has disappeared
\begin{equation}\label{LimitField}
u_{0}\in \mH^{1}_{0}(\Omega)
\quad\textrm{ and }\quad
-\Delta u_{0} = f\quad\textrm{in}\;\Omega .
\end{equation}
However, in the general case (more precisely, when $u_0(0)\ne0$), the convergence turns out to be very 
slow: for any arbitrary radius $r>0$ such that $\overline{\mD}_r\subset\Om$, we have 
$\|u_{\delta}-u_{0}\|_{\mH^{1}(\Om\setminus\overline{\mD}_r)}\geq\,C\,|\ln\delta|^{-1}\|f\|_{\mrm{L}^2(\Om)}$, for some constant 
$C>0$ independent of $\delta$. To give an idea  $|\ln\delta|^{-1}\approx 0.0434$ for $\delta=10^{-10}$. 
Thus, neglecting the presence of the small inclusion is not satisfactory from a computational point of view,
and a reasonable numerical approach for (\ref{PbInit}) should reproduce accurately the perturbation 
induced by the presence of the small inclusion $\omega_{\delta}$. \\
\newline
Most of the numerical approaches that could be considered for dealing with this problem
suffer from a numerical \textit{locking effect} \cite{BaSu92}: performances of standard strategies deteriorate as $\delta\to 0$.
Admittedly,  robust strategies already exist in the literature, like the multi-scale finite element 
method coupled with some mesh refinement strategy \cite{MR2477579}, or the boundary element method (see 
e.g. \cite{MR3058835}). These techniques provide satisfying results in many cases, but they require careful and 
thorough implementation efforts, and/or rely on strong assumptions such as homogeneity of the coefficients 
of the equation under consideration. Other numerical strategies are based on an approximation of $u_{\delta}$ 
of the form "$u_{0}$ + \textrm{corrector}", where both terms of this sum are computed separately (see for example 
\cite{DaVi05,BDTV07}). These approaches may induce substantial additional computational cost in a real life simulation. 
In many practical situations, small inclusions are not the main subject of concern, and it would be desirable 
to devise a simple, general purpose and implementation friendly method that would not rely on any kind of mesh 
refinement technique, while remaining robust as $\delta\to 0$. This is the purpose of the present article to 
describe and analyse a method matching these requirements, while relying on only one numerical resolution.\\
\newline
The outline of this article is the following. In Section \ref{sectionAsyExp}, we summarize the main results concerning 
the asymptotic expansion of $u_{\delta}$ with respect to the size of the small hole.  Section \ref{sectionConstruction} 
is dedicated to the construction of a model problem, based on the matched expansion of $u_{\delta}$, whose solution 
has the same far field asymptotics as $u_{\delta}$, up to a remainder in $O(\delta^{1-\epsilon}),\forall\epsilon>0$. 
We prove this in Section \ref{sectionAnalysis} (see Proposition  \ref{EstimateApprox}) and also show that 
consistency of any Galerkin discretization of this model problem is quasi-optimal and uniform with respect
to $\delta$. Finally, in Section \ref{sectionNum}, we present and comment numerical results that confirm and 
illustrate our theoretical conclusions.

\section{Asymptotic expansion of the solution}\label{sectionAsyExp}
Asymptotic analysis for problems involving small inclusions can be found in many works. We refer the reader 
to \cite{MR1429996,MR630701,Ilin92,Na61,MR912054,MaNP00,Naza99,NaSo03,NaSo06}. The asymptotic expansion for the particular problem 
we are considering in this paper is described in detail in \cite{MaNP00} and here, we just wish to remind the 
main results provided by the method of matched expansions at order one. To proceed, we need first 
to introduce two particular functions: the Green function $G$ and the logarithmic capacity potential $P$. These 
two functions are defined in normalized geometries by the following equations 
\begin{equation}\label{Profiles}
\left\{\begin{array}{ll}
-\Delta G  =  0 & \textrm{in}\;\;\Omega\setminus\{O\}\\[5pt]
G = 0 & \textrm{on}\;\;\partial\Omega\\
\multicolumn{2}{l}{\dsp{
G(\bfx) = \frac{1}{2\pi}\ln(1/ \vert\bfx\vert)+
\mathop{O}_{\vert\bfx\vert\to 0}(1)}}
\end{array}\right.
\quad\quad\quad
\left\{\begin{array}{ll}
-\Delta P = 0 &\textrm{in}\;\;\Xi:=\R^{2}\setminus\overline{\omega}\\[5pt]
P = 0&\textrm{on}\;\;\partial\Xi\\
\multicolumn{2}{l}{\dsp{ P(\bfxi)=  \frac{1}{2\pi}\ln(1/ \vert\bfxi\vert)+
\mathop{O}_{\vert\bfxi\vert\to \infty}(1)}.}
\end{array}\right.
\end{equation}
Classical techniques of separation of variables (see e.g. \cite{MR0350075}) show that there 
exist constants $G_{0}$, $P_{0}$ that depend only on the domains $\Omega$, $\omega$ such that 
$G(\bfx) -(2\pi)^{-1}\ln\vert\bfx\vert^{-1} - G_{0}  = O(\vert \bfx\vert)$ for 
$\vert \bfx\vert\to 0$, and $P(\bfxi) -(2\pi)^{-1}\ln\vert\bfxi\vert^{-1} - P_{0}  = O(\vert \bfxi\vert^{-1})$
for $\vert\bfxi\vert\to\infty $. The asymptotic analysis of Problem (\ref{PbInit}) also involves 
the gauge function (see \cite{MR630701})
\begin{equation}\label{deflambda}
\lambda(\delta) := \dsp\frac{2\pi}{\ln\delta+2\pi(P_{0}-G_{0})  }\ .
\end{equation}
Finally, the global approximation of $u_{\delta}$ is defined as an interpolation 
between a far field and a near field contribution as follows:
\begin{equation}\label{MatchedAsymptotics1}
\hat{u}_{\delta}(\bfx)\  :=\ \psi(\bfx/\delta)
\,v_{\delta}(\bfx) + \chi(\bfx)\,V_{\delta}(\bfx/\delta) 
-\chi(\bfx)\psi(\bfx/\delta)\,m_{\delta}(\bfx) 
\end{equation}
\begin{equation}\label{MatchedAsymptotics}
\begin{array}{l|lcl}
\mbox{where}\  ~& v_{\delta}(\bfx) & := & u_{0}(\bfx) +u_{0}(0)\,
\lambda(\delta)\,G(\bfx) \\[10pt]
&V_{\delta}(\bfxi) & := & u_{0}(0)\,\lambda(\delta)\,P(\bfxi)  \\[6pt]
&m_{\delta}(\bfx) & := & \dsp u_{0}(0)\,\lambda(\delta)\,\Big( 
\frac{1}{2\pi}\ln(\delta/\vert\bfx\vert) + P_{0} \Big).\phantom{spaacee}
\end{array}
\end{equation}
In the expression above, the cut-off functions $\chi$, $\psi$ are two elements of $\mathscr{C}^{\infty}(\overline{\Om}):= \{v\vert_{\Omega}\,\vert\, v\in\mathscr{C}^{\infty}(\R^2) \}$ 
such that $\chi(\bfx) = 1$ for $\vert\bfx\vert \leq r_{0}/2$, $\chi(\bfx) = 0$ for $\vert\bfx\vert \geq r_{0}$ and $\psi:= 
1 - \chi$. Here, $r_0>0$ is a given parameter such that $\mrm{D}_{r_0}\subset \Omega$. The following well-known result 
provides an error estimate for $\Vert u_{\delta} - \hat{u}_{\delta}\Vert_{\mH^{1}(\Omega_{\delta})}$. For the proof, we refer the 
reader, for example, to Section 2.4.1 of \cite{MaNP00}. 
\begin{proposition}\quad\\
Considering $u_{\delta}$ defined by (\ref{PbInit}) and $\hat{u}_{\delta}$ defined by (\ref{MatchedAsymptotics}), there exist 
constants $C$, $\delta_0>0$ independent of $\delta$ such that 
$$
\Vert u_{\delta}-\hat{u}_{\delta} \Vert_{\mH^{1}(\Omega_{\delta})}\leq C\,\delta\,
\vert\ln\delta\vert\,\|f\|_{\mrm{L}^2(\Om)}\qquad\forall \delta\in (0,\delta_0].
$$
\end{proposition}
\noindent  Note that, the constant $C$ involved in the estimate above 
a priori depends on $\chi,\psi$. Note also that, in the definition of $\chi$, $\psi$, 
the parameter $r_{0}$ could be \textit{any} positive number such that 
$\overline{\mrm{D}}_{r_0}\subset\Om$. In particular, it can be chosen arbitrarily small. 
Looking at the explicit definition of $\hat{u}_{\delta}$ given by (\ref{MatchedAsymptotics}), 
this implies the following result.
\begin{proposition}\label{FarFieldEstimate}\quad\\
Consider $u_{\delta}$, $v_{\delta}$ defined by (\ref{PbInit}), (\ref{MatchedAsymptotics}). For any disk 
$\mD_r\subset \Omega$ with $0<r\le r_0$, there exist constants $C_{r}$, $\delta_0>0$ independent of $\delta$ such that 
\begin{equation}\label{approxFF}
\Vert\, u_{\delta}-v_{\delta}\,\Vert_{\mH^{1}(\Omega\setminus\overline{\mD}_r)}\leq 
C_{r}\,\delta\,\vert\ln\delta\vert\,\|f\|_{\mrm{L}^2(\Om)}\qquad\forall \delta\in (0,\delta_0].
\end{equation}
\end{proposition}

\noindent This last result shows that $v_{\delta}=u_{0} + u_{0}(0)\lambda(\delta)G$ provides a reasonable 
approximation (for example $\delta\,\vert\ln\delta\vert\approx 2.3\ 10^{-9}$ for $\delta=10^{-10}$)
of $u_{\delta}$ at any fixed distance from the small hole. Thus, the far field of $u_{\delta}$ appears as the 
superposition of the limit field $u_{0}$ and a field ``radiated'' by a point source located at the center of the hole. \\
Numerically, $u_{0}$, $G$ can be approximated by functions $u_{0}^{h}$, $G^{h}$ using a standard finite element 
method (here, $h$ refers to some mesh size) and define $v_{\delta}^{h}=u_{0}^{h} + u_{0}^{h}(0)\lambda(\delta)G^{h}$. Then, 
(\ref{approxFF}) ensures that $v_{\delta}^{h}$ is a good approximation of the far field of $u_{\delta}$. This procedure 
is rather simple to implement and it has been proven in \cite{BoDa13} (see also \cite{BDTV09,DaVi05,DaVi07,BDTV07,BHBo05,BoVo00} 
for slightly different problems\footnote{This technique is also very close to singular complement methods (or singular function methods) which are used to compute efficiently the solution of elliptic partial differential equations in non 
smooth domains (see for example \cite{BDLN92,CiLa11,CJKLZ05a,CJKLZ05b,HaLo02})}) that it gives good results. However, it requires to solve two problems which we 
would like to avoid because it may be time consuming. Adapting this approach to the case of $N$ inclusions 
would lead to $N+1$ numerical solves. Similarly, looking for an approximation of $u_{\delta}$ as sharp as the 
first $M$ terms of its asymptotic expansion would lead to $M$ numerical solves. From this perspective, for 
practical computations, a method involving only one numerical solve would be much more interesting. \\ 
In the next section, we propose a model problem that can be discretized by means of any standard 
Galerkin method (with classical finite elements for example) with quasi-optimal approximation properties
with respect to both $h$ and $\delta$. In addition, the numerical schemes obtained in this manner 
do not deteriorate as  $\delta\to 0$. 

\section{Construction of a model problem}\label{sectionConstruction}
The model problem we wish to propose is formulated in (\ref{VarFor}). The goal of the present section is to explain how we obtain 
this problem. To avoid having to compute both $u_{0}$ and $G$ in the decomposition $v_{\delta}=u_{0} + u_{0}(0)\lambda(\delta)G$, 
we will use the fact that the regular part of $G$ belongs to $\mH^1(\Om)$. Let us decompose $G$ under the form
\begin{equation}\label{NewDecompG}
G=\bfs_{\log} + \tilde{G}\qquad\mbox{with}\qquad \bfs_{\log}(\bfx):=\frac{1}{2\pi}\,\chi(\bfx)\ln(1/\vert\bfx\vert)\quad
\mbox{ and }\quad\tilde{G}\in\mH^{1}_{0}(\Omega)\cap\mathscr{C}^{0}(\Omega).
\end{equation}
This allows us to write $v_{\delta}$ as $v_{\delta} = w_{\delta}+ u_{0}(0)\lambda(\delta)\,\bfs_{\log}$ with $w_{\delta}:=u_{0} + 
u_{0}(0)\lambda(\delta)\tilde{G}$. Let us express the coefficient $u_{0}(0)\lambda(\delta)$ by means of $w_{\delta}$. 
According to the definition of $u_{0}$ and $G$ it is clear that $w_{\delta}$ belongs 
to $\mH^{1}_{0}(\Omega)\cap\mathscr{C}^{0}(\Om)$, where $\mathscr{C}^{0}(\Om)$ refers to the space of continuous functions
on $\overline{\Omega}$. Moreover, observing that $\tilde{G}=G_0+\hat{G}$ for some function 
$\hat{G}\in \mH^{1}_{0}(\Omega)\cap\mathscr{C}^{0}(\Om)$ 
vanishing at $0$, we find $w_{\delta}(0)=u_{0}(0)(1+\lambda(\delta)\,G_0)$ and so 
$u_{0}(0)\lambda(\delta)=w_{\delta}(0)\lambda(\delta)/(1+\lambda(\delta)\,G_0)$. Using Definition (\ref{deflambda}) of $\lambda(\delta)$, 
we deduce that 
\begin{equation}\label{decompoBis}
v_{\delta} = w_{\delta}+ b_{\delta}(w_{\delta})\,\bfs_{\log}\qquad\textrm{with}\qquad w_{\delta}:=u_{0} + u_{0}(0)\lambda(\delta)\tilde{G}
\quad\mbox{ and } \quad b_{\delta}(w_{\delta}): =  \dsp\frac{2\pi\,w_{\delta}(0)}{\ln\delta + 2\pi P_{0}}.
\end{equation}
Let us emphasize that this expression for $v_{\delta}$ is interesting because it involves only one unknown function which belongs 
to the variational space $\mH^{1}_{0}(\Omega)$. Now, we need to derive a problem characterizing $w_{\delta}$. In the 
sense of distributions in $\Omega$, there holds $-\Delta w_{\delta}=-\Delta (u_{0} + u_{0}(0)\lambda(\delta)\,\tilde{G})=
f-b_{\delta}(w_{\delta})\,\Delta \tilde{G}$. Multiplying by $w'\in\mH^{1}_{0}(\Omega)$ and using Green's formula, we find that $w_{\delta}$ verifies
\begin{equation}\label{VarId}
\begin{array}{l}
a(w_{\delta},w')+b_{\delta}(w_{\delta})\,b_{\log}(w')  = \dsp\int_{\Omega}
fw'\,d\bfx 
\end{array}
\end{equation}
\begin{equation}\label{DefTerVaria}
\begin{array}{l|lcl}
\mbox{where}\qquad  ~& a(w_{\delta},w') & := & \dsp \int_{\Omega}\nabla w_{\delta}\cdot\nabla w' \,d\bfx \\[10pt]
& b_{\log}(w') & := & \dsp\int_{\Omega} \Delta \tilde{G}\,w' \,d\bfx  \\[6pt]
&  \Delta \tilde{G}(\bfx) & := &  (2\pi)^{-1}\Big((\Delta\chi)(\bfx)\ln\vert\bfx\vert
 + 2\nabla\chi(\bfx)\cdot\nabla(\ln\vert\bfx\vert ) \Big).\phantom{spaacee}
\end{array}
\end{equation}
Note that $\chi$ is equal to one in a neighbourhood of $0$ so that $\Delta \tilde{G}$ indeed belongs to $\mathscr{C}^{\infty}(\overline{\Om})$. As a remark, let us observe that for test functions $w'$ such that 
$0\notin \mrm{supp}(w')$, we have $b_{\log}(w') = \int_{\Omega}\nabla 
\bfs_{\log}\cdot\nabla w'\,d\bfx$.  \\
\newline
Of course (\ref{VarId}) is not a valid variational formulation in $\mH^{1}(\Omega)$ because the functional 
$w_{\delta}\mapsto w_{\delta}(0)$ is not defined on this space. So it cannot be exploited directly for discretization and 
then numerical computation. This is the motivation for considering a regularized counterpart of (\ref{VarId}) 
where in $b_{\delta}(w_{\delta})$, we replace $w_{\delta}(0)$ by $(2\pi\delta)^{-1}\int_{\partial\mD_{\delta}}w_{\delta}\,d\sigma$, 
$\partial\mD_{\delta}$ denoting the circle centered at $0$ and of radius $\delta$. Finally, this leads us 
to examine the following model problem, 
\begin{equation}\label{VarFor}
\begin{array}{|l}
\textrm{Find $\tilde{w}_{\delta}\in\mH^{1}_{0}(\Omega)$ such that}\\
a(\tilde{w}_{\delta},w') + \tilde{b}_{\delta}(\tilde{w}_{\delta})\,b_{\log}(w') =
\dsp\int_{\Omega}f w'\,d\bfx\qquad \forall w'\in\mH^{1}_{0}(\Omega),  
\end{array}
\end{equation}
where $a(\cdot,\cdot)$, $b_{\log}(\cdot)$ are defined in (\ref{DefTerVaria}) and where 
\begin{equation}\label{def_b_delta}
\tilde{b}_{\delta}(\tilde{w}_{\delta})\;:=\;\dsp{\frac{2\pi}{\ln\delta + 
2\pi P_{0}}}\,\frac{1}{2\pi\delta}\int_{\partial\mD_{\delta}}\tilde{w}_{\delta}\,d\sigma.
\end{equation}
The variational formulation (\ref{VarFor}) perfectly makes sense for $\delta$ small enough and, 
in the next  section, we show that it admits a unique solution so that  $\tilde{w}_{\delta}$ is well defined.
Since Problem (\ref{VarFor}) differs from (\ref{VarId}), its solution $\tilde{w}_{\delta}$ is a priori different from $w_{\delta}$. 
However we are going to show that $\tilde{w}_{\delta}$ and $w_{\delta}$ (defined by (\ref{decompoBis})) are close to each other, and that  
$\tilde{w}_{\delta} + \tilde{b}_{\delta}(\tilde{w}_{\delta})\,\bfs_{\log}$ is a good approximation of the far field of $u_{\delta}$. 

\begin{remark}
We could have proposed a formulation where, in (\ref{VarId}), the term $w_{\delta}(0)$ is replaced  by $(\pi\delta^2)^{-1}\int_{\mD_{\delta}}w_{\delta}\,d\bfx$. The analysis we will develop and the results we will obtain would have been the same with this alternative choice. 

\end{remark}
\begin{remark}
It is worth noting that in (\ref{VarFor}), a simple perturbation of a usual formulation allows to take into account the small hole. Therefore, 
with this approach, we can adapt classical codes at little cost. In this respect, this technique shares similarities with the extended finite element method 
(XFEM) \cite{BeBl99,DoBe99} and the generalized finite element method (GFEM) \cite{Duar96,MeBa96}. 
\end{remark}
\begin{remark}
In this paper, we have chosen to investigate the problem of the small hole with Dirichlet boundary condition in 2D only because in this case, 
the logarithmic term which appears in the asymptotic expansion of $u_{\delta}$ makes the zero order approximation clearly unsatisfactory (see 
the discussion in the introduction). However, the present approach could allow to consider other 
problems of singular perturbation. It could also be adapted to obtain higher orders of approximation. In this case, new perturbation terms 
would have to be considered in the left hand side of (\ref{VarFor}).
\end{remark}

\section{Analysis and discretization of the model problem}\label{sectionAnalysis}

We first prove that $\tilde{a}_{\delta}(\cdot,\cdot):=a(\cdot,\cdot)+\tilde{b}_{\delta}(\cdot)\,b_{\log}(\cdot)$, the bilinear form appearing in 
the left hand side of (\ref{VarFor}), differs from  $a(\cdot,\cdot)$ by a small perturbation. This will allow to show that $\tilde{w}_{\delta}$ is a 
relevant approximation of $w_{\delta}$.
\begin{proposition}\label{Perturbation}\quad\\
There exists a constant $C>0$ independent of $\delta$ such that 
\begin{equation}\label{UniformContinuity}
\sup_{\varphi\in \mH^{1}_{0}(\Omega)\setminus\{0\}}
\frac{\vert \tilde{b}_{\delta}(\varphi)\vert}{\Vert \varphi\Vert_{\mH^{1}(\Omega)}}
\leq \frac{C}{\sqrt{\vert\ln\delta\vert}}\quad\quad\forall\delta\in (0,1). 
\end{equation}
As a consequence, for $\delta$ small enough, $\tilde{a}_{\delta}(\cdot,\cdot)$ is coercive and Problem (\ref{VarFor}) has a unique solution 
$\tilde{w}_{\delta}$. Moreover, for any $\eps>0$, there exist constants $C_{\eps}$, $\delta_0>0$  independent of $\delta$ such that
\begin{equation}\label{EstimateW}
\|\tilde{w}_{\delta}-w_{\delta}\|_{\mH^{1}(\Omega)}\leq \,C_{\eps}\,\delta^{1-\eps}\,\|f\|_{\mrm{L}^2(\Om)}\qquad\forall\delta\in(0,\delta_0] ,
\end{equation}
where $w_{\delta}$ is the function defined in (\ref{decompoBis}).
\end{proposition}
\noindent \textbf{Proof:} First, we prove (\ref{UniformContinuity}). Consider the disk $\mD_{r_0}$  introduced in the 
definition of $\chi$ that satisfies  $\overline{\mD}_{r_0}\subset \Omega$. We have in particular 
$\supp(\chi)\subset \mD_{r_0}$. Take an arbitrary $\zeta\in \mathscr{C}^{\infty}(\overline{\Omega})$ such 
that $\supp (\zeta)\subset \overline{\mD}_{r_0}$. Integration by parts and Cauchy-Buniakowski-Schwarz 
inequality show that
\begin{equation}\label{estimHardy}
\dsp\Big\vert\frac{1}{2\pi\delta}\int_{\partial\mD_{\delta}} \zeta\,d\sigma \Big\vert  =
\Big\vert\frac{1}{2\pi}\int_{\mD_{r_0}\setminus\overline{\mD}_{\delta}} \nabla(\ln \vert\bfx\vert)\cdot\nabla\zeta\,d\bfx \Big\vert\ 
\leq\ \sqrt{\frac{\vert \ln(r_{0}/\delta) \vert}{2\pi}}\;\Vert \zeta\Vert_{\mH^{1}(\Omega)} \;.
\end{equation}
As a consequence, for any $\varphi\in \mathscr{C}^{\infty}_{0}(\Omega)$, considering $\chi\varphi$ instead of $\zeta$ in (\ref{estimHardy}) and 
using (\ref{def_b_delta}), we see that there exist constants $C,C'>0$ (whose values may change from one occurrence to another) independent of $\delta$ such that 
$\vert \tilde{b}_{\delta}(\varphi)\vert = \vert \tilde{b}_{\delta}(\chi\varphi)\vert\leq C \vert \ln\delta\vert^{-1/2}
\,\Vert \chi\varphi\Vert_{\mH^{1}(\Omega)}$ $\leq C' \vert \ln\delta\vert^{-1/2}\,\Vert \varphi\Vert_{\mH^{1}(\Omega)}$. Since $\mathscr{C}^{\infty}_{0}(\Omega)$ is dense into $\mH^{1}_{0}(\Omega)$, this shows (\ref{UniformContinuity}). We deduce that for all $\varphi\in\mH^{1}_0(\Omega)$, we have
\begin{equation}\label{proofCoer}
|\tilde{a}_{\delta}(\varphi,\varphi)| = |a(\varphi,\varphi)+\tilde{b}_{\delta}(\varphi)\,b_{\log}(\varphi)| \ge C\,(1-C'\,|\ln \delta|^{-1/2}) \|\varphi \|^2_{\mH^{1}(\Omega)}.
\end{equation}
This guarantees that for $\delta$ small enough, Problem (\ref{VarFor}) has a unique solution $\tilde{w}_{\delta}$.
To establish the second part of the statement, we use (\ref{proofCoer}) and write, for $\delta$ small enough, 
\begin{equation}\label{equationInter1}
\dsp\Vert  w_{\delta} - \tilde{w}_{\delta}\Vert^2_{\mH^{1}(\Omega)}\leq C\,|\tilde{a}_{\delta}(w_{\delta} - \tilde{w}_{\delta},w_{\delta} - \tilde{w}_{\delta})|\le C\,|b_{\delta}(w_{\delta})-\tilde{b}_{\delta}(w_{\delta})|\  |b_{\log}(w_{\delta} - \tilde{w}_{\delta})|.
\end{equation}
Let us focus on
\begin{equation}\label{equationDiff}
|b_{\delta}(w_{\delta})-\tilde{b}_{\delta}(w_{\delta})|=\Big\vert\dsp\frac{2\pi}{\ln\delta+2\pi\,P_0}\Big\vert\,\Big\vert\, w_{\delta}(0)-\frac{1}{2\pi\delta}\,\int_{\partial\mD_{\delta}}w_{\delta}\,d\sigma\,\Big\vert.
\end{equation}
From (\ref{decompoBis}), we know that there holds $w_{\delta} = u_{0} + u_{0}(0)\lambda(\delta)\tilde{G}$ with $\tilde{G}=G_0+\hat{G}$, $\hat{G}\in\mathscr{C}^{\infty}(\overline{\Om})$. 
For $\beta\in \R$, we define the weighted norm  
\begin{equation}\label{WeightedNorm}
\Vert \varphi\Vert_{\mV^{1}_{\beta}(\Omega)} := (\ \Vert \,\vert\bfx\vert^\beta \nabla\varphi\Vert^2_{\mrm{L}^2(\Omega)}+\Vert \,\vert\bfx\vert^{\beta-1}\varphi\Vert^2_{\mrm{L}^2(\Omega)}\ )^{1/2},
\end{equation}
and let $\mV^{1}_{\beta}(\Omega)$ refer to the completion of $\mathscr{C}^{\infty}(\overline{\Omega}\setminus\{O\})$
$:= \{ v\vert_{\Omega}\;\vert\; v\in\mathscr{C}^{\infty}(\R^{2}),\,v = 0$ in a neighbourhood of $0$ $\}$ 
with respect to this norm. We refer the reader to \cite{Kond67} for more details on weighted Sobolev spaces.
On the other hand, classical Kondratiev analysis (see \cite[Chap.6]{KoMR97}) allows to prove the decomposition 
$u_{0} = u_{0}(0)+\tilde{u}_{0}$, where $\tilde{u}_{0}\in \mH^{1}(\Omega)\cap \mV^{1}_{-1+\eps}(\Omega)$ 
for all $\eps>0$, with the estimate
\begin{equation}\label{estimKondra}
|u_{0}(0)| + \|\tilde{u}_{0}\|_{\mV^{1}_{-1+\eps}(\Omega)}\le C_{\eps}\,\|f\|_{\mrm{L}^2(\Om)}.
\end{equation}
This implies $w_{\delta} = w_{\delta}(0)+(\tilde{u}_{0} + u_{0}(0)\lambda(\delta)\hat{G})$. Conducting a calculus analogue to (\ref{estimHardy}), 
replacing formally $\zeta(\bfx)$ by $\vert\bfx\vert^{-1+\eps}\tilde{u}_{0}(\bfx)$, we find that $|\int_{\partial\mD_{\delta}} \tilde{u}_{0}\,d\sigma| \le  C\,\delta^{2-\eps}\,\|f\|_{\mrm{L}^2(\Om)}$. 
Writing a Taylor expansion of $\hat{G}$ at $\bfx=0$ and using (\ref{estimKondra}), we obtain $|u_{0}(0)\lambda(\delta)\int_{\partial\mD_{\delta}} 
\hat{G}\,d\sigma| \le  C\,\delta^{2}\,\|f\|_{\mrm{L}^2(\Om)}$. We deduce 
\begin{equation}\label{estimPoids}
\Big\vert\, w_{\delta}(0)-\frac{1}{2\pi\delta}\,\int_{\partial\mD_{\delta}}w_{\delta}\,d\sigma\,\Big\vert \le C\,\delta^{1-\eps}\,\|f\|_{\mrm{L}^2(\Om)}.
\end{equation}
Plugging this estimate in (\ref{equationDiff}) leads to $|b_{\delta}(w_{\delta})-\tilde{b}_{\delta}(w_{\delta})| \leq \,C\,\delta^{1-\eps}\,\|f\|_{\mrm{L}^2(\Om)}$.
Combining this inequality with (\ref{equationInter1}), we obtain  (\ref{EstimateW}) as a direct consequence.\hfill $\Box$\\

\noindent We have just proved that $\tilde{w}_{\delta}$ is close to $w_{\delta}$. From the relation linking $w_{\delta}$ to 
$v_{\delta}$, we deduce that $\tilde{w}_{\delta} + \tilde{b}_{\delta}(\tilde{w}_{\delta})\,\bfs_{\log}$ is a good approximation of the far field of $u_{\delta}$.
\begin{proposition}\label{EstimateApprox}\quad\\
For any disk $\mD_r\subset \Omega$ with $0<r\le r_0$ and for any $\eps>0$, there exist constants $C$, $\delta_0>0$ depending on $r$, $\eps$ 
but not on $\delta$ such that
\begin{equation}\label{approxFFW}
\Vert u_{\delta} - (\tilde{w}_{\delta} + \tilde{b}_{\delta}(\tilde{w}_{\delta})\,\bfs_{\log})
\Vert_{\mH^{1}(\Omega\setminus\overline{\mD}_r)}\leq C\,\delta^{1-\eps}\,\|f\|_{\mrm{L}^2(\Om)}\qquad\forall \delta\in(0,\delta_0].
\end{equation}
\end{proposition}
\noindent \textbf{Proof:} Remembering that $v_{\delta} = w_{\delta}+ b_{\delta}(w_{\delta})$ (see (\ref{decompoBis})), 
where $v_{\delta}$, $w_{\delta}$ are defined in (\ref{MatchedAsymptotics}), (\ref{decompoBis}), and using the triangular inequality, we can write
\begin{equation}\label{ErrorEstimate1}
\begin{array}{l}
\phantom{\le\ }  \Vert u_{\delta} - (\tilde{w}_{\delta} + \tilde{b}_{\delta}(\tilde{w}_{\delta})\,\bfs_{\log})
\Vert_{\mH^{1}(\Omega\setminus\overline{\mD}_r)} \\[5pt] 
 \le  \Vert u_{\delta} - v_{\delta}\Vert_{\mH^{1}(\Omega\setminus\overline{\mD}_r)}+\Vert w_{\delta} - \tilde{w}_{\delta}\Vert_{\mH^{1}(\Omega\setminus\overline{\mD}_r)}+
|b_{\delta}(w_{\delta})-\tilde{b}_{\delta}(\tilde{w}_{\delta})|\,\Vert \bfs_{\log}
\Vert_{\mH^{1}(\Omega\setminus\overline{\mD}_r)}.
\end{array}
\end{equation}
In the previous proof, we have established that  $|b_{\delta}(w_{\delta})-\tilde{b}_{\delta}(w_{\delta})| \leq \,C\,\delta^{1-\eps}\,\|f\|_{\mrm{L}^2(\Om)}$ 
for some constant $C>0$ independent of $\delta$. Combining this with (\ref{UniformContinuity}), we find
\begin{equation}\label{equationInter2}
\begin{array}{lcl}
|b_{\delta}(w_{\delta})-\tilde{b}_{\delta}(\tilde{w}_{\delta})| & \le & |b_{\delta}(w_{\delta})-\tilde{b}_{\delta}(w_{\delta})|+|\tilde{b}_{\delta}(w_{\delta})-\tilde{b}_{\delta}(\tilde{w}_{\delta})| \\[5pt]
& \le & C\,(\delta^{1-\eps}\,\|f\|_{\mrm{L}^2(\Om)}+\Vert w_{\delta} - \tilde{w}_{\delta}\Vert_{\mH^{1}(\Omega)}),
\end{array}
\end{equation}
for some constant $C>0$ independent of $\delta$. Plugging (\ref{equationInter2}) in (\ref{ErrorEstimate1}) and using 
(\ref{approxFF}), (\ref{EstimateW}), we finally obtain (\ref{approxFFW}).\hfill $\Box$
\begin{remark}\quad\\
Working as in the previous proof, one can obtain a slightly more general result where the norm $\|\cdot\|_{\mH^{1}(\Omega\setminus\overline{\mD}_r)}$ 
in the right hand side of (\ref{approxFFW}) is replaced by the norm $\|\cdot\|_{\mV^{1}_{\beta}(\Omega)}$ (see (\ref{WeightedNorm})) with $\beta>0$.
\end{remark}
\begin{remark}\label{RmqRegu}\quad\\
Making the additional assumption that the source term $f$ verifies $\|\,\vert\bfx\vert^{-\beta}f\|_{\mrm{L}^2(\Om)}<+\infty$ for some $\beta>0$, 
and revisiting Estimate (\ref{estimPoids}), we find that (\ref{EstimateW}) can be improved in 
$\|\tilde{w}_{\delta}-w_{\delta}\|_{\mH^{1}(\Omega)}\leq \,C\,\delta\,\|\,\vert\bfx\vert^{-\beta}f\|_{\mrm{L}^2(\Om)}$, $\forall\delta\in(0,\delta_0]$. 
In this case, (\ref{approxFFW}) becomes $\Vert u_{\delta} - (\tilde{w}_{\delta} + \tilde{b}_{\delta}(\tilde{w}_{\delta})\,\bfs_{\log})
\Vert_{\mH^{1}(\Omega\setminus\overline{\mD}_r)}\leq C\,\delta\,\vert\ln\delta\vert\,\|\,\vert\bfx\vert^{-\beta}f\|_{\mrm{L}^2(\Om)}$, $\forall\delta\in(0,\delta_0]$.
\end{remark}

\noindent To conclude, assume that we want to solve Formulation (\ref{VarFor}) by means of 
a Galerkin approach associated with a family of discrete subspaces $(\mV^{h})_{h>0}$ (in the numerical experiments, $h$ will refer to the mesh size). We assume that there holds $\mV^{h}\subset 
\mH^{1}_{0}(\Omega)$ for all $h>0$ . The natural discrete variational formulation associated with (\ref{VarFor}) writes
\begin{equation}\label{DiscreteVF}
\begin{array}{|l}
\mbox{Find }\tilde{w}_{\delta}^{h}\in \mV^{h}\mbox{ such that}\\
a(\tilde{w}_{\delta}^{h},\varphi^{h}) + \tilde{b}_{\delta}(\tilde{w}_{\delta}^{h})\,b_{\log}(\varphi^{h}) = 
\dsp\int_{\Omega}f \varphi^{h}\,d\bfx\qquad\forall \varphi^{h}\in \mV^{h}.
\end{array}
\end{equation}
The coercivity of $\tilde{a}_{\delta}(\cdot,\cdot)=a(\cdot,\cdot)+\tilde{b}_{\delta}(\cdot)\,b_{\log}(\cdot)$ proven in 
Proposition \ref{Perturbation} shows straightforwardly, by Cea's lemma, the result of quasi-optimal convergence  
\begin{equation}\label{quasiOptimalCv}
\Vert \tilde{w}_{\delta}-\tilde{w}_{\delta}^{h}\Vert_{\mH^{1}(\Omega)}\leq C
\inf_{\varphi^{h}\in \mV^{h}}\Vert \tilde{w}_{\delta}-\varphi^{h}\Vert_{\mH^{1}(\Omega)}.
\end{equation}
Combining this with Estimate (\ref{approxFFW}) proves that $\tilde{w}_{\delta}^{h}+\tilde{b}_{\delta}(\tilde{w}_{\delta}^{h})\,\,\bfs_{\log}$ 
is a reasonable approximation of the far field expansion of $u_{\delta}$. The following proposition is one of the 
two main results (with Proposition \ref{PropoMainBis} hereafter) of the present article. It establishes quasi-optimal convergence of the numerical method (\ref{DiscreteVF})
\textit{both} in $\delta$ and $h$.

\begin{proposition}\label{Corol}\quad\\
Consider a finite dimensional space $\mV^{h}\subset \mH^{1}_{0}(\Omega)$. For any disk $\mD_r\subset \Omega$ with 
$0<r\le r_0$ and for any $\eps>0$, there exists a constant $C>0$ depending on $r$, $\eps$ but not on $\delta$ 
and $h$ such that, for $\delta$ small enough,
\begin{equation}\label{QuasiUniformConv}
\begin{array}{l}
  \Vert u_{\delta} -( \tilde{w}_{\delta}^{h} + \tilde{b}_{\delta}(\tilde{w}_{\delta}^{h})\,\bfs_{\log} )
  \Vert_{\mH^{1}(\Omega\setminus\overline{\mD}_r)}\\[10pt]
  \quad\quad \leq  \dsp{  C\,( \delta^{1-\eps}+|\ln\delta|^{-1}\,\mathop{\inf}_{\varphi^{h}\in \mV^{h}}\Vert \tilde{G}-\varphi^{h}
    \Vert_{\mH^{1}(\Omega)})\,\|f\|_{\mrm{L}^2(\Om)}+ C \mathop{\inf}_{\varphi^{h}\in \mV^{h}}\Vert u_{0}-\varphi^{h}\Vert_{\mH^{1}(\Omega)} }.
\end{array}
\end{equation}
\end{proposition}
\noindent \textbf{Proof:} 
The continuity estimate of $\tilde{b}_{\delta}$ (see (\ref{UniformContinuity})) implies that there exists a constant 
$C>0$ independent of $\delta$ such that  $\Vert u_{\delta} -(\tilde{w}_{\delta}^{h}+\tilde{b}_{\delta}(\tilde{w}_{\delta}^{h}) 
\bfs_{\log})\Vert_{\mH^{1}(\Omega\setminus\overline{\mD}_r)}\leq C\,\Vert u_{\delta} - (\tilde{w}_{\delta}+\tilde{b}_{\delta}(\tilde{w}_{\delta}) 
\bfs_{\log})\Vert_{\mH^{1}(\Omega\setminus\overline{\mD}_r)}$   
$ + C\,\Vert\tilde{w}_{\delta}-\tilde{w}^h_{\delta}\Vert_{\mH^{1}(\Omega)}$. Proposition \ref{EstimateApprox} already 
yields that $\Vert u_{\delta} - (\tilde{w}_{\delta} + \tilde{b}_{\delta}(\tilde{w}_{\delta})\,\bfs_{\log})
\Vert_{\mH^{1}(\Omega\setminus\overline{\mD}_r)}\leq C\,\delta^{1-\eps}\,\|f\|_{\mrm{L}^2(\Om)}$ for any $\eps>0$, so we 
only need to focus on the second term of the previous inequality. Since we have 
$w_{\delta} = u_{0}+u_{0}(0)\lambda(\delta) \tilde{G}$, (\ref{quasiOptimalCv}) allows us to write
$$
\begin{array}{ll}
 & \Vert\tilde{w}_{\delta}-\tilde{w}_{\delta}^{h}\Vert_{\mH^{1} (\Omega)}\\[5pt]
\leq & C\,\Vert w_{\delta}- \tilde{w}_{\delta}\Vert_{\mH^{1}(\Omega)} + 
C \inf_{\varphi^{h}\in \mV^{h}}\Vert w_{\delta}-\varphi^{h}\Vert_{\mH^{1}(\Omega)} \\[5pt]

\leq &  
C\,\Vert w_{\delta}- \tilde{w}_{\delta}\Vert_{\mH^{1}(\Omega)} + 
C \inf_{\varphi^{h}\in \mV^{h}}\Vert u_{0}-\varphi^{h}\Vert_{\mH^{1}(\Omega)} + 
C\,\vert u_{0}(0)\lambda(\delta)\vert\,\inf_{\varphi^{h}\in \mV^{h}}\Vert \tilde{G}-\varphi^{h}
\Vert_{\mH^{1}(\Omega)} \\[5pt]

\leq &  C\,( \delta^{1-\eps}+|\ln\delta|^{-1}\,\inf_{\varphi^{h}\in \mV^{h}}\Vert \tilde{G}-\varphi^{h}
\Vert_{\mH^{1}(\Omega)})\,\|f\|_{\mrm{L}^2(\Om)}+ 
C \inf_{\varphi^{h}\in \mV^{h}}\Vert u_{0}-\varphi^{h}\Vert_{\mH^{1}(\Omega)}.
\end{array}
$$
This finishes the proof.
\hfill $\Box$\\
\newline
To illustrate what kind of result the above proposition implies, assume for example that $\mV^{h}$ is the space of 
$\mathbb{P}_{1}$-Lagrange finite element functions constructed on a quasi-uniform regular triangulation of the 
domain $\Omega$. In this situation, according to (\ref{QuasiUniformConv}), for any $\eps >0$ there exists a 
constant $C>0$ independent of $\delta$ and $h$, such that $\Vert u_{\delta} -( \tilde{w}_{\delta}^{h} + 
\tilde{b}_{\delta}(\tilde{w}_{\delta}^{h})\,\bfs_{\log} )\Vert_{\mH^{1}(\Omega\setminus\overline{\mD}_r)}\leq C\, (\delta^{1-\eps}+h)\,\|f\|_{\mrm{L}^2(\Om)}$.

We also emphasize that the result of quasi-optimal convergence for (\ref{DiscreteVF}) with a constant independent of 
$\delta$ discards any numerical locking effect. In other words, this assert the robustness of (\ref{DiscreteVF}) as 
$\delta\to 0$.

\section{Practical implementation of the perturbation}

From the point of view of practical implementation, a natural idea consists in computing the perturbation 
term $\tilde{b}_{\delta}(\cdot)$ by means of the crude quadrature formula $\int_{\partial\mrm{D}_{\delta}}\varphi_{h}d\sigma 
\simeq 2\pi\delta\varphi_{h}(0)$ for any $\varphi_{h}\in\mV_{h}$, which boils down to actually considering 
$b_{\delta}(\cdot)$ instead of $\tilde{b}_{\delta}(\cdot)$. In this section we examine the validity of such a 
substitution. We introduce the discrete formulation 
\begin{equation}\label{DiscreteVFBis}
\begin{array}{|l}
\mbox{Find }w_{\delta}^{h}\in \mV^{h}\mbox{ such that}\\
a(w_{\delta}^{h},\varphi^{h}) + b_{\delta}(w_{\delta}^{h})\,b_{\log}(\varphi^{h}) = 
\dsp\int_{\Omega}f \varphi^{h}\,d\bfx\qquad\forall \varphi^{h}\in \mV^{h}.
\end{array}
\end{equation}
assuming that $\mV^{h}\subset\mathscr{C}^0(\overline{\Om})$ (this implies in particular that (\ref{DiscreteVFBis}) has 
indeed a  sense) is a Lagrange finite element space constructed on a quasi-uniform regular triangulation of the 
domain $\Omega$. Let us prove that Problem (\ref{DiscreteVFBis}) yields to a good approximation of the far field 
of $u_{\delta}$. 
\begin{proposition}\label{PropoMainBis}\quad\\
For any given $h>0$, for $\delta>0$ small enough, Problem (\ref{DiscreteVFBis}) has a unique solution $w_{\delta}^{h}$. Moreover, if $f\in\mH^2(\Om)$ and if $\Om$ is smooth, then for any disk $\mD_r\subset \Omega$ with 
$0<r\le r_0$ and for any $\eps>0$, there exists a constant $C>0$ depending on $r$, $\eps$ but not on $\delta$ 
and $h$ such that, for $\delta$ small enough,
\begin{equation}\label{QuasiUniformConvBis}
\begin{array}{ll}
  \quad \Vert u_{\delta} -( w_{\delta}^{h} + b_{\delta}(w_{\delta}^{h})\,\bfs_{\log} )
  \Vert_{\mH^{1}(\Omega\setminus\overline{\mD}_r)}\\[10pt]
  \leq   \dsp C\,( \delta|\ln\delta|+\gamma(\delta,h)+|\ln\delta|^{-1}\,\mathop{\inf}_{\varphi^{h}\in \mV^{h}}\Vert \tilde{G}-\varphi^{h}
    \Vert_{\mH^{1}(\Omega)})\,\|f\|_{\mH^2(\Om)}+ C \mathop{\inf}_{\varphi^{h}\in \mV^{h}}\Vert u_{0}-\varphi^{h}\Vert_{\mH^{1}(\Omega)} .
\end{array}
\end{equation}
In (\ref{QuasiUniformConvBis}), the constant $\gamma(\delta,h)$ can be chosen such that $\gamma(\delta,h)=(\delta+h^2|\ln h|)/(1-(1+|\ln h|)^{1/2}/|\ln\delta|)$. 
\end{proposition}
\begin{remark}
Observe that for any given $h>0$, there holds $|\gamma(\delta,h)|\le C\,(\delta+h^2|\ln h|)$ for $\delta$ small enough. Actually, in the proof, we will see that the condition $|\ln h|^{1/2}/|\ln\delta|=O(1)$ is sufficient to guarantee well-posedness for Problem (\ref{DiscreteVFBis}). Note that this assumption is in accordance with the situation we want to consider, namely an obstacle small compare to the mesh size ($\delta<<h$).
\end{remark}
\begin{remark}
The additional smoothness assumption on the source term is needed for technical reasons (see the proof of Lemma \ref{lemmaUniConv}). The authors do not know if it can be weakened.
\end{remark}

\noindent \textbf{Proof:} We first recall the discrete Sobolev inequality (see \cite[Lemma 4.9.2]{BrSc08})
\begin{equation}\label{DiscreteSoboIneq}
\|\varphi^{h}\|_{\mrm{L}^{\infty}(\Om)} \le C\,(1+|\ln h|)^{1/2})\,\|\varphi^{h}\|_{\mH^1(\Om)}\qquad \forall \varphi^{h}\in \mV^{h}.
\end{equation}
Here and in the sequel of this proof, $C>0$ denotes a constant independent of $\delta, h$ which may change from one occurrence to another. Since $b_{\delta}(\varphi^{h}) = 2\pi\,\varphi^{h}(0)/(\ln\delta + 2\pi P_{0})$, we deduce from (\ref{DiscreteSoboIneq}) that, for $\delta$ small enough, for all $\varphi^{h}\in \mV^{h}$, we have
\begin{equation}\label{estimCoer}
|a(\varphi^{h},\varphi^{h})+b_{\delta}(\varphi^{h})\,b_{\log}(\varphi^{h})| \ge C\,\alpha(\delta,h)\,\|\varphi^{h}\|^2_{\mH^1(\Om)},
\end{equation}
where $\alpha(\delta,h):=1-\beta(\delta,h)$ and $\beta(\delta,h):=(1+|\ln h|)^{1/2}/|\ln\delta|$. It is clear that for a given $h$, $\beta(\delta,h)$ tends to zero as $\delta$ goes to zero. Therefore, Estimate (\ref{estimCoer}) shows that $a(\cdot,\cdot)+b_{\delta}(\cdot)\,b_{\log}(\cdot)$ is coercive for $\delta$ small enough. In this case, from the Lax-Milgram theorem, we infer that Problem (\ref{DiscreteVFBis}) has a unique solution. Now, we wish to establish (\ref{QuasiUniformConvBis}). Thanks to the triangular inequality, we can write
\begin{equation}\label{termInter1}
\begin{array}{l}
  \Vert u_{\delta} -( w_{\delta}^{h} + b_{\delta}(w_{\delta}^{h})\,\bfs_{\log} )
  \Vert_{\mH^{1}(\Omega\setminus\overline{\mD}_r)} \le \Vert u_{\delta} -( \tilde{w}_{\delta}^{h} + \tilde{b}_{\delta}(\tilde{w}_{\delta}^{h})\,\bfs_{\log} )
  \Vert_{\mH^{1}(\Omega\setminus\overline{\mD}_r)}\\[4pt]
\hspace{5.4cm}+\|w_{\delta}^{h}-\tilde{w}_{\delta}^{h}\|_{\mH^{1}(\Omega\setminus\overline{\mD}_r)}+|b_{\delta}(w_{\delta}^h)-\tilde{b}_{\delta}(\tilde{w}_{\delta}^h)|\,\|\bfs_{\log}\|_{\mH^{1}(\Omega\setminus\overline{\mD}_r)}.
\end{array}
\end{equation}
The first term of the right hand side of (\ref{termInter1}) has already been studied in Proposition \ref{Corol}. To handle the last term, we use (\ref{DiscreteSoboIneq}) to obtain
\begin{equation}\label{termInter2}
\begin{array}{lcl}
|b_{\delta}(w_{\delta}^h)-\tilde{b}_{\delta}(\tilde{w}_{\delta}^h)| & \le & |b_{\delta}(w_{\delta}^h)-b_{\delta}(\tilde{w}_{\delta}^h)|+|b_{\delta}(\tilde{w}_{\delta}^h)-\tilde{b}_{\delta}(\tilde{w}_{\delta}^h)|\\[4pt]
& \le & C\,\beta(\delta,h)\,\|w_{\delta}^{h}-\tilde{w}_{\delta}^{h}\|_{\mH^{1}(\Omega)}+|b_{\delta}(\tilde{w}_{\delta}^h)-\tilde{b}_{\delta}(\tilde{w}_{\delta}^h)|.
\end{array}
\end{equation}
Let us estimate the quantity $\|w_{\delta}^{h}-\tilde{w}_{\delta}^{h}\|_{\mH^{1}(\Omega)}$ which appears both in (\ref{termInter1}) and (\ref{termInter2}). The coercivity inequality (\ref{estimCoer}) and the definition of Problems (\ref{DiscreteVF}), (\ref{DiscreteVFBis}) provide 
\[
\begin{array}{lcl}
C\,\alpha(\delta,h)\,\|w_{\delta}^{h}-\tilde{w}_{\delta}^{h}\|_{\mH^{1}(\Omega)}^2  & \le & |a(w_{\delta}^{h}-\tilde{w}_{\delta}^{h},w_{\delta}^{h}-\tilde{w}_{\delta}^{h})+b_{\delta}(w_{\delta}^{h}-\tilde{w}_{\delta}^{h})\,b_{\log}(w_{\delta}^{h}-\tilde{w}_{\delta}^{h})| \\[4pt]
& \le & |b_{\delta}(\tilde{w}_{\delta}^h)-\tilde{b}_{\delta}(\tilde{w}_{\delta}^h)|\,|b_{\log}(\tilde{w}_{\delta}^{h}-\tilde{w}_{\delta}^{h})|.
\end{array}
\]
Observing that $b_{\log}$ is bounded on $\mH^1(\Om)$, we deduce that 
\begin{equation}\label{termInter3}
\|w_{\delta}^{h}-\tilde{w}_{\delta}^{h}\|_{\mH^{1}(\Omega)}\le C\,\alpha(\delta,h)^{-1}\,|b_{\delta}(\tilde{w}_{\delta}^h)-\tilde{b}_{\delta}(\tilde{w}_{\delta}^h)|.
\end{equation}
Plugging (\ref{termInter3}) in (\ref{termInter1}) and (\ref{termInter2}), we conclude that it is sufficient to control $|b_{\delta}(\tilde{w}_{\delta}^h)-\tilde{b}_{\delta}(\tilde{w}_{\delta}^h)|$ to prove (\ref{QuasiUniformConvBis}).  We have
\begin{equation}\label{estimdiffb}
|b_{\delta}(\tilde{w}_{\delta}^h)-\tilde{b}_{\delta}(\tilde{w}_{\delta}^h)| \le |b_{\delta}(\tilde{w}_{\delta}-\tilde{w}_{\delta}^h)-\tilde{b}_{\delta}(\tilde{w}_{\delta}-\tilde{w}_{\delta}^h)|+|b_{\delta}(\tilde{w}_{\delta})-\tilde{b}_{\delta}(\tilde{w}_{\delta})|.
\end{equation}
Then, by definition of $b_{\delta}$, $\tilde{b}_{\delta}$, we find, for $\delta$ small enough,
\begin{equation}\label{estimDiffForm}
\begin{array}{lcl}
|b_{\delta}(\tilde{w}_{\delta}-\tilde{w}_{\delta}^h)-\tilde{b}_{\delta}(\tilde{w}_{\delta}-\tilde{w}_{\delta}^h)| & = & \Big|\dsp{\frac{2\pi}{\ln\delta + 
2\pi P_{0}}}\,\frac{1}{2\pi\delta}\int_{\partial\mD_{\delta}}(\tilde{w}_{\delta}-\tilde{w}_{\delta}^h)(0)-(\tilde{w}_{\delta}-\tilde{w}_{\delta}^h)
\,d\sigma\Big|\\[12pt]
& \le & C\,\|\tilde{w}_{\delta}-\tilde{w}_{\delta}^h\|_{\mrm{L}^{\infty}(\Om)}
\end{array}
\end{equation}
Lemma \ref{lemmaUniConv} hereafter guarantees that if $\tilde{w}_{\delta}\in \mrm{W}^{2,\infty}(\Om):=\{v\in\mrm{L}^{\infty}(\Om)\,|\,\partial^{\boldsymbol{\alpha}}v\in\mrm{L}^{\infty}(\Om),\,|\boldsymbol{\alpha}|\le2\}$\footnote{In this definition, we use the classical \textit{multi-index} notation.}, then $(\tilde{w}_{\delta}^h)$ uniformly converges to $\tilde{w}_{\delta}$ as $h$ tends to zero, with the estimate 

\begin{equation}\label{EstimUnifConv}
\|\tilde{w}_{\delta}-\tilde{w}_{\delta}^h\|_{\mrm{L}^{\infty}(\Om)}\le C\,h^2 |\ln h|\,\|\tilde{w}_{\delta}\|_{\mrm{W}^{2,\infty}(\Om)}.
\end{equation}
To ensure such a regularity for $\tilde{w}_{\delta}$, let us assume that the source term $f$ verifies $f\in\mH^2(\Om)$. Using (\ref{VarFor}) and (\ref{def_b_delta}), we see that in the sense of distributions in $\Omega$, there holds
\[
-\Delta \tilde{w}_{\delta}=f_{\delta}\qquad\mbox{with}\qquad f_{\delta}:=f-\tilde{b}_{\delta}(\tilde{w}_{\delta})\Delta \tilde{G}.
\]
Thus, if $f\in\mH^2(\Om)$, then the theory of elliptic regularity asserts that $\tilde{w}_{\delta}\in\mH^4(\Om)$. Besides, Proposition \ref{Perturbation} and Estimate (\ref{proofCoer}) imply
\begin{equation}\label{eqEllipticReg}
\|\tilde{w}_{\delta}\|_{\mH^1(\Om)}\le C\,\|f\|_{\mrm{L}^2(\Om)}.
\end{equation}
We emphasize that in (\ref{eqEllipticReg}), the constant $C>0$ is independent of $\delta$. As a consequence, from Proposition \ref{Perturbation} and (\ref{eqEllipticReg}), we get $\|f_{\delta}\|_{\mH^2(\Om)} \le C\,\|f\|_{\mH^2(\Om)}$. In this case, from the Sobolev imbedding theorem, we deduce that $\tilde{w}_{\delta}\in\mathcal{C}^2(\overline{\Om})$ with  
\begin{equation}\label{eqEllipticReg2}
\|\tilde{w}_{\delta}\|_{\mrm{W}^{2,\infty}(\Om)}\le C\,\|\tilde{w}_{\delta}\|_{\mH^4(\Om)} \le C\,\|f_{\delta}\|_{\mH^2(\Om)} \le C\,\|f\|_{\mH^2(\Om)}.
\end{equation}
Collecting (\ref{estimDiffForm}), (\ref{EstimUnifConv}) and (\ref{eqEllipticReg2}), we find 
\begin{equation}\label{eqEllipticReg3}
|b_{\delta}(\tilde{w}_{\delta}-\tilde{w}_{\delta}^h)-\tilde{b}_{\delta}(\tilde{w}_{\delta}-\tilde{w}_{\delta}^h)| \le C\,h^2 |\ln h|\,\|f\|_{\mH^2(\Om)}.
\end{equation}
On the other hand, concerning the second term of the right hand side of (\ref{estimdiffb}), writing the Taylor expansion of $\tilde{w}_{\delta}$ at $\bfx=0$ and coming back to the definition of $b_{\delta}$, $\tilde{b}_{\delta}$ yields 
\begin{equation}\label{eqEllipticReg4}
|b_{\delta}(\tilde{w}_{\delta})-\tilde{b}_{\delta}(\tilde{w}_{\delta})| \le C\,\delta\,\|f\|_{\mH^2(\Om)}.
\end{equation}
Plugging (\ref{eqEllipticReg3}) and  (\ref{eqEllipticReg4}) into (\ref{estimdiffb}) leads to 
\begin{equation}\label{FinalEstimate}
|b_{\delta}(\tilde{w}_{\delta}^h)-\tilde{b}_{\delta}(\tilde{w}_{\delta}^h)| \le C\,(\delta+h^2 |\ln h|)\,\|f\|_{\mH^2(\Om)}.
\end{equation}
Finally, combining (\ref{termInter1}), (\ref{termInter2}), (\ref{termInter3}), (\ref{FinalEstimate}) and using Remark \ref{RmqRegu} allows to obtain (\ref{QuasiUniformConvBis}). \hfill $\Box$\\
\newline
In order to complete the previous analysis, now we state a result of uniform approximation of $\tilde{w}_{\delta}$ by $\tilde{w}_{\delta}^h$ whose proof can be obtained working exactly as in \cite{Scot76}.

\begin{lemma}\label{lemmaUniConv}
Assume that the solution of Problem (\ref{VarFor}) verifies $\tilde{w}_{\delta}\in\mrm{W}^{2,\infty}(\Om)$. Then, for $\delta$ small enough, we have the estimate
\[
\|\tilde{w}_{\delta}-\tilde{w}_{\delta}^h\|_{\mrm{L}^{\infty}(\Om)}\le C\,h^2 |\ln h|\,\|\tilde{w}_{\delta}\|_{\mrm{W}^{2,\infty}(\Om)},
\]
where $C>0$ is independent of $\delta,h$. 
\end{lemma}

\section{Numerical experiments}\label{sectionNum}
\begin{figure}[!ht]
\centering
\includegraphics[scale=0.65]{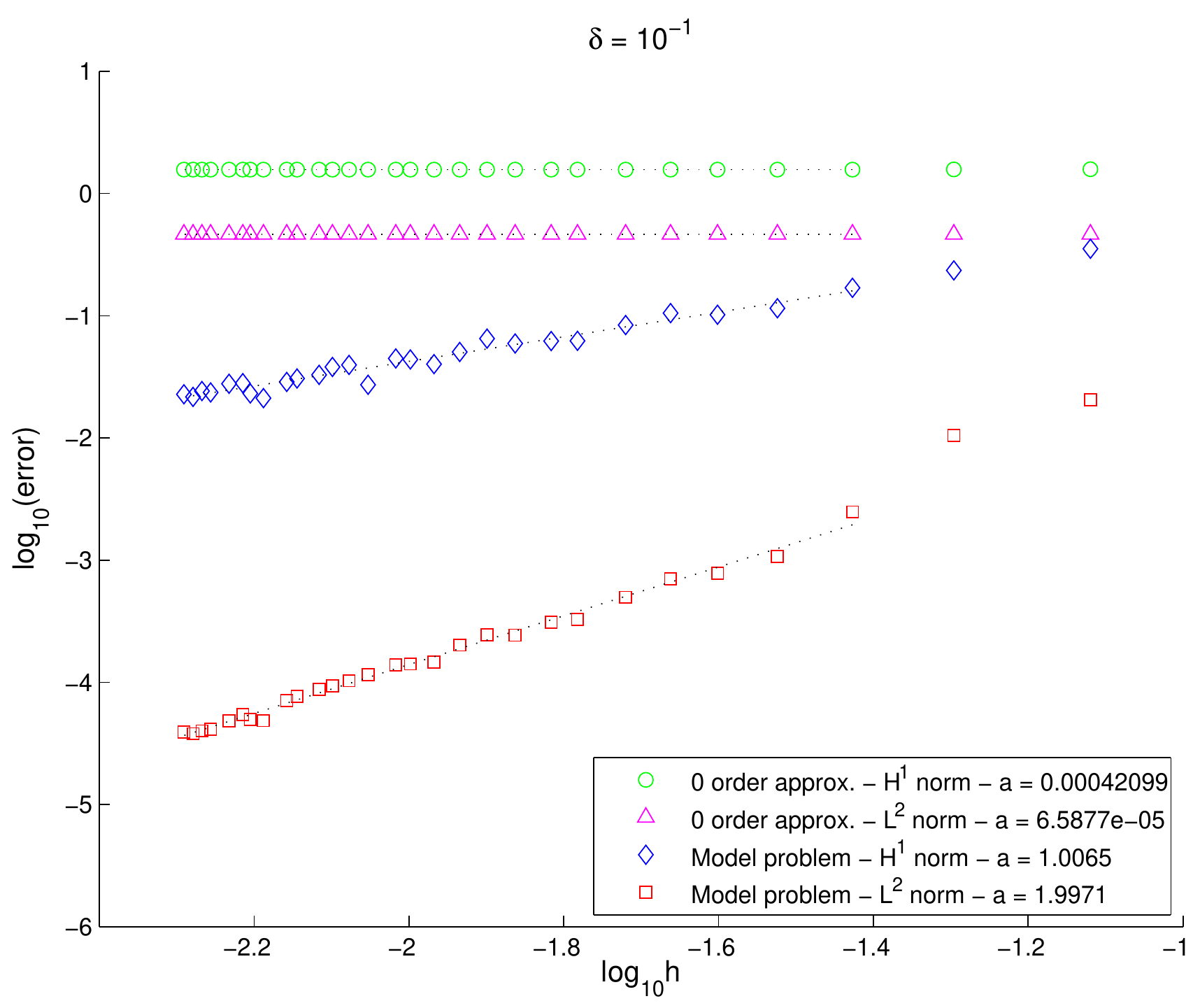}
\caption{Convergence w.r.t. the mesh size -- $\delta=10^{-1}$, $g=1$. \label{Results1}}
\end{figure}

\begin{figure}[!ht]
\centering
\includegraphics[scale=0.65]{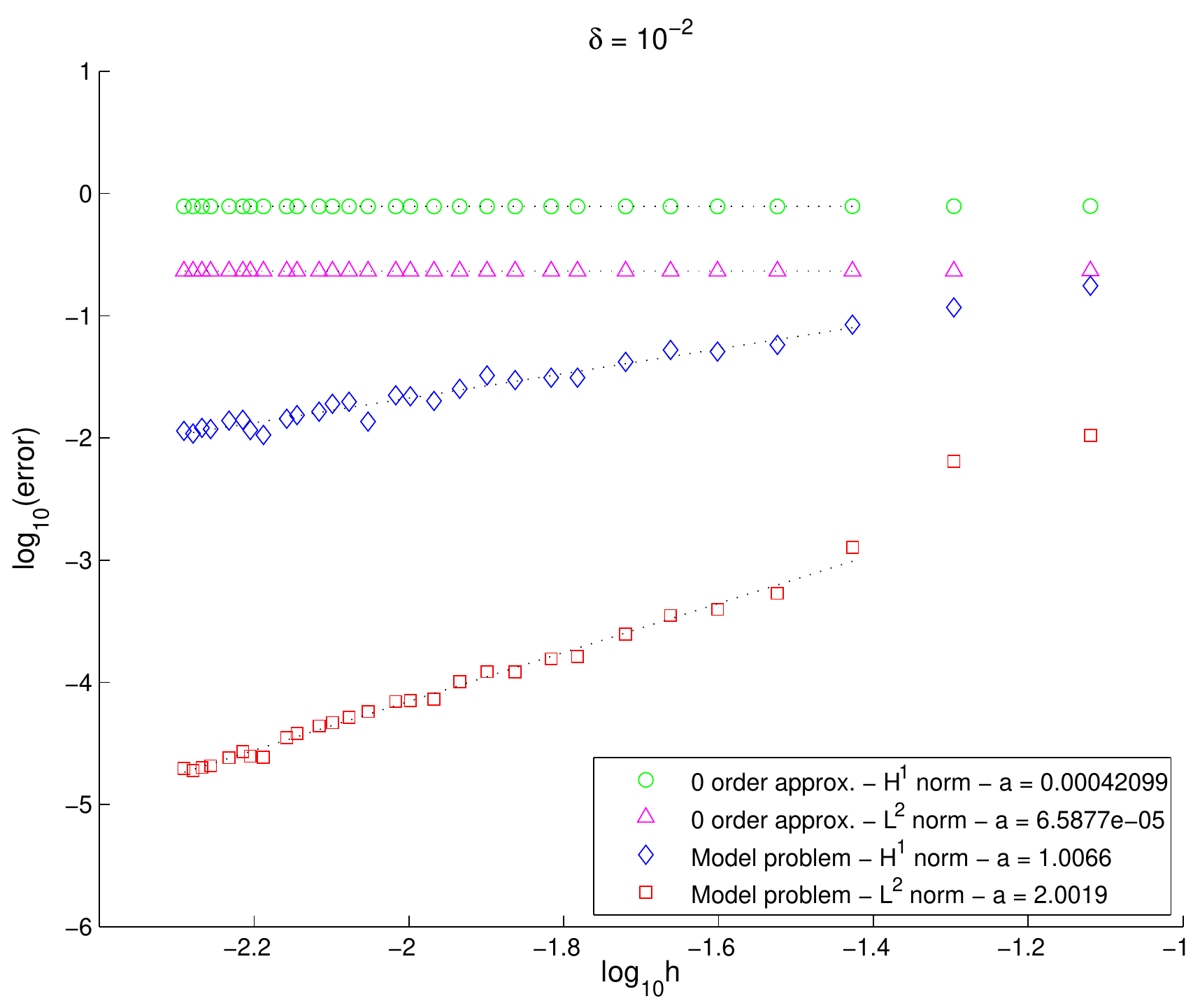}
\caption{Convergence w.r.t. the mesh size -- $\delta=10^{-2}$, $g=1$.\label{Results2}}
\end{figure}

\begin{figure}[!ht]
\centering
\includegraphics[scale=0.65]{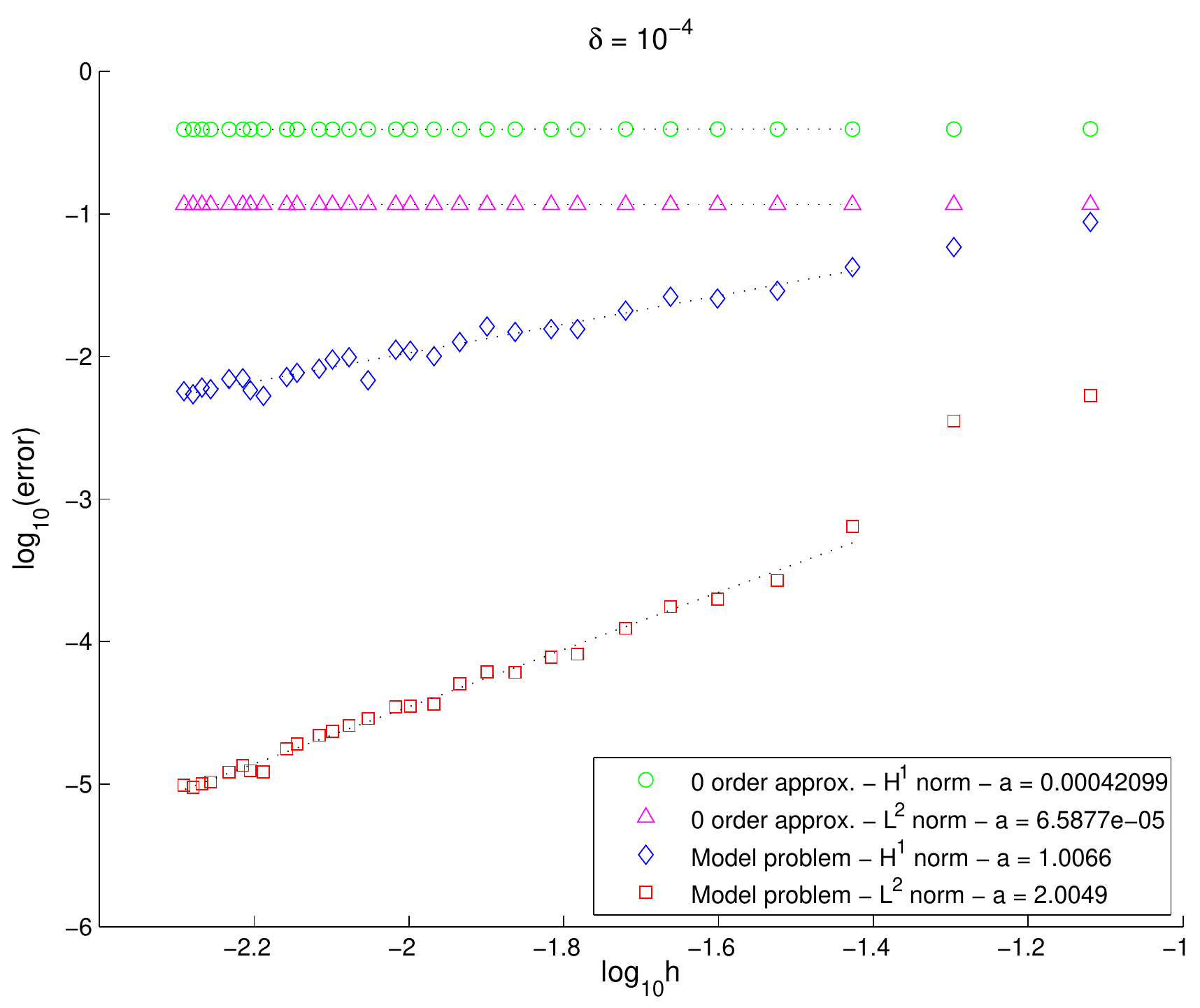}
\caption{Convergence w.r.t. the mesh size -- $\delta=10^{-4}$, $g=1$.\label{Results3}}
\end{figure}

\begin{figure}[!ht]
\centering
\includegraphics[scale=0.65]{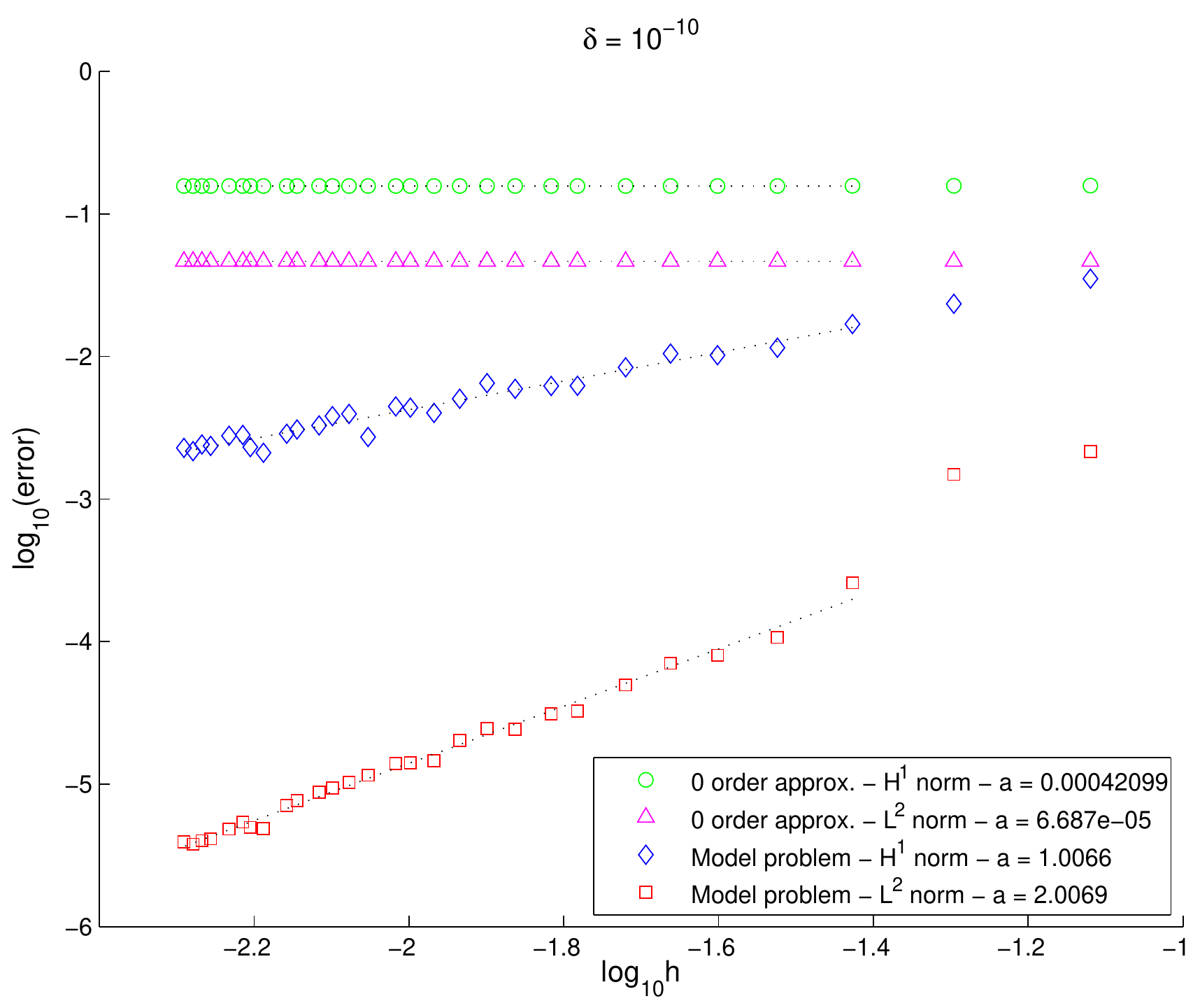}
\caption{Convergence w.r.t. the mesh size -- $\delta=10^{-10}$, $g=1$.\label{Results4}}
\end{figure}

\begin{figure}[!ht]
\centering
\includegraphics[scale=0.65]{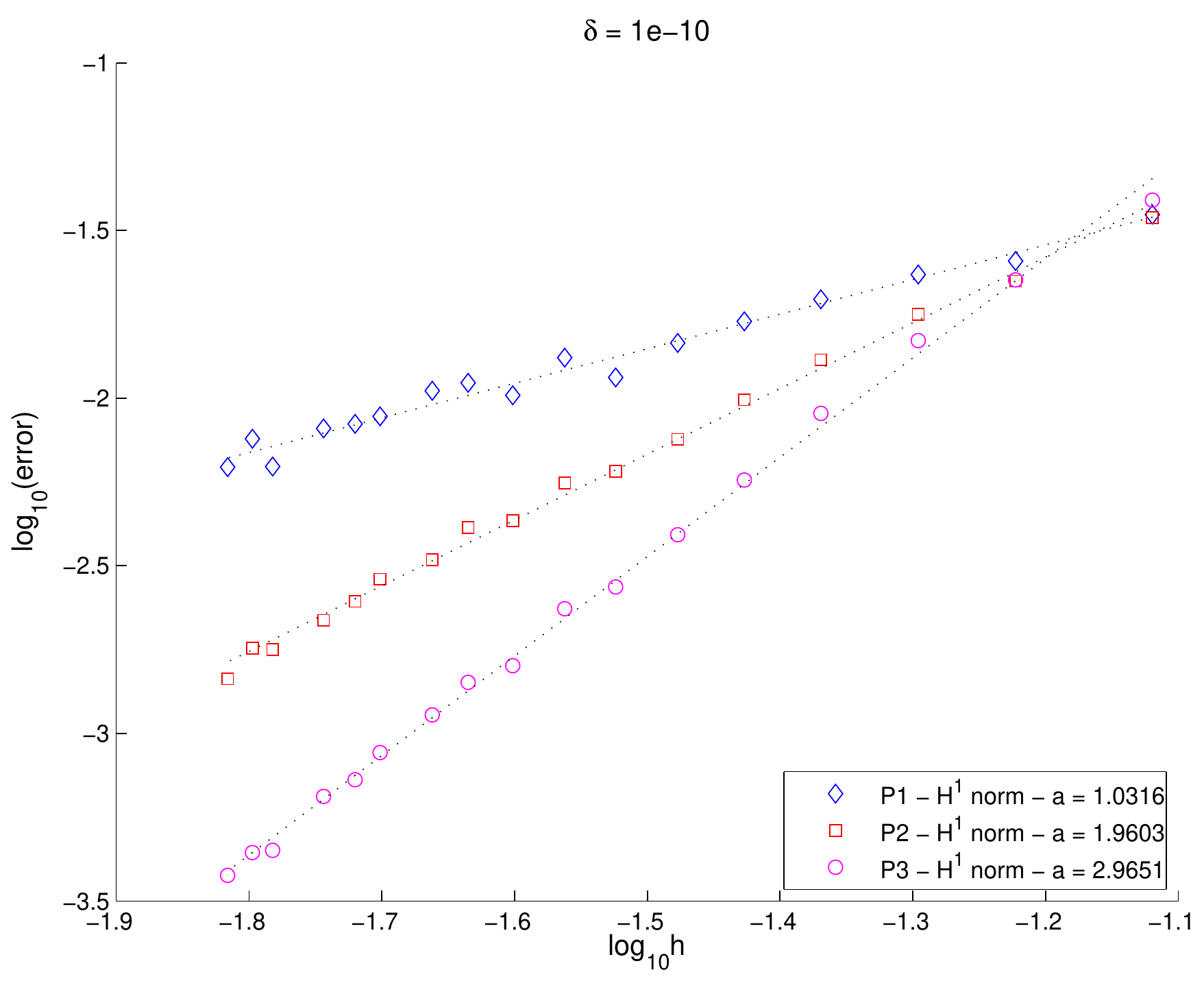}
\caption{Convergence w.r.t. the mesh size for several orders of approximation -- $\delta=10^{-10}$, $g=1$. The cut-off function $\chi$ appearing in the definition of $b_{\log}(\cdot)$ (see (\ref{DefTerVaria})) is chosen in  $\mathscr{C}^{\infty}(\overline{\Om})$.\label{Orders}}
\end{figure}

\begin{figure}[!ht]
\centering
\includegraphics[scale=0.65]{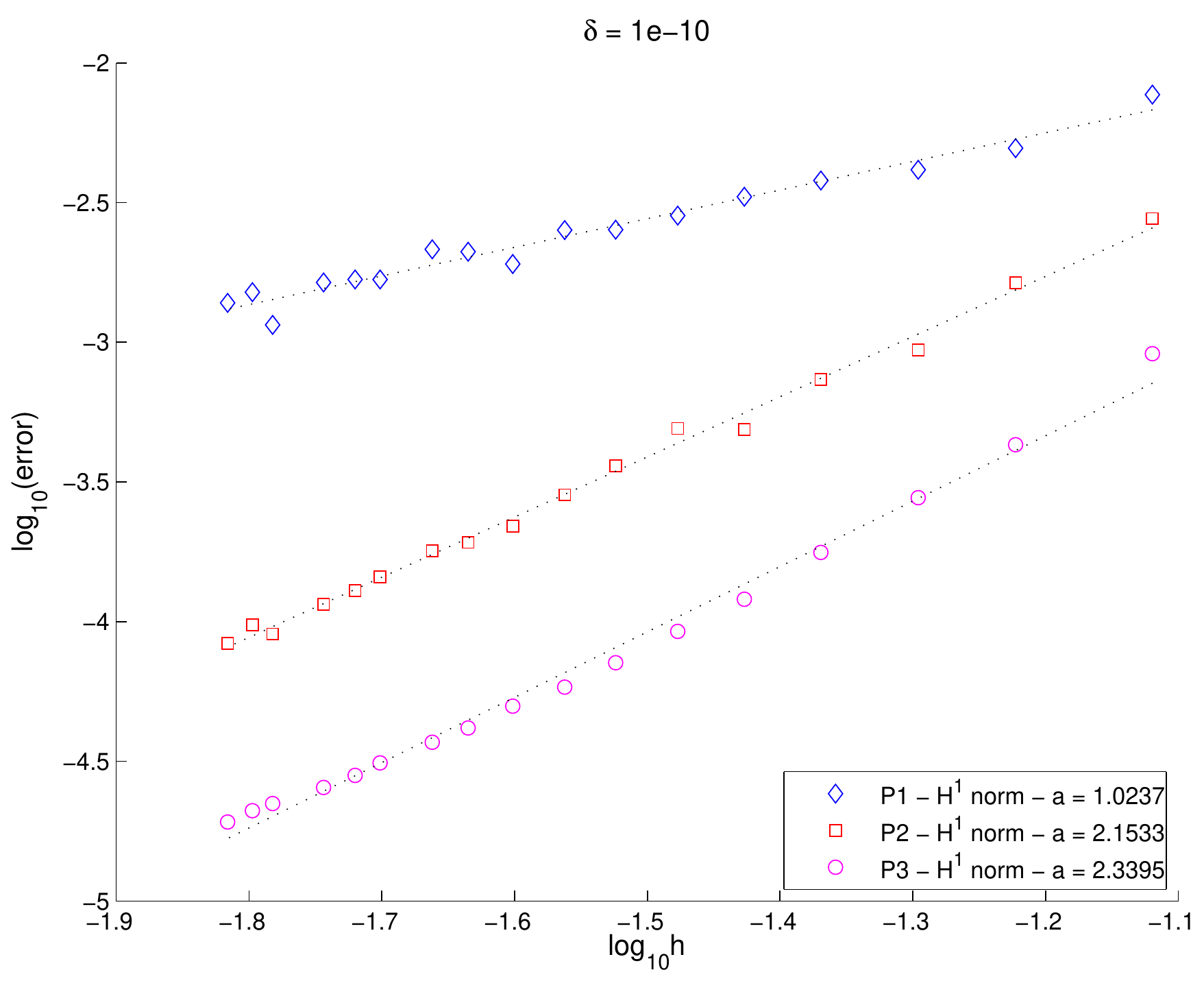}
\caption{Convergence w.r.t. the mesh size for several orders of approximation -- $\delta=10^{-10}$, $g=1$. The cut-off function $\chi$ appearing in the definition of $b_{\log}(\cdot)$ (see (\ref{DefTerVaria})) is chosen in  $\mathscr{C}^{3}(\overline{\Om})\setminus\mathscr{C}^{4}(\overline{\Om})$.\label{OrdersC3}}
\end{figure}

\begin{figure}[!ht]
\centering
\includegraphics[scale=0.65]{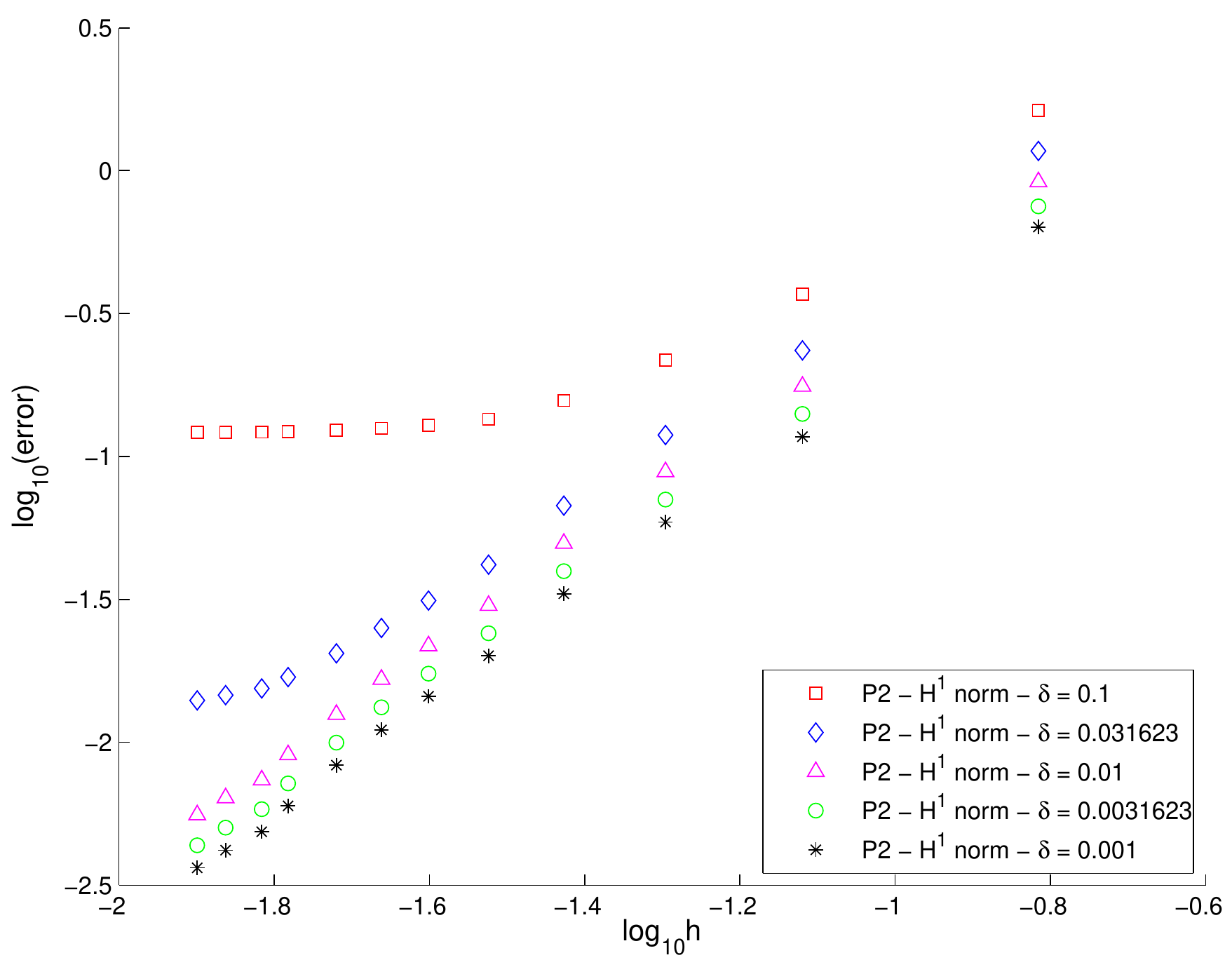}
\caption{Convergence w.r.t. the mesh size for several values of $\delta$ -- $g=1+\sin(\theta)$.\label{ResultsThreshold}}
\end{figure}

Now, let us present the results of the numerical tests that we conducted in order to validate our theoretical conclusions. First,
we detail the parameters used for the experiments. Let $\Om$ (resp. $\om_{\delta}$) be the disk centered at $0$ of radius $1$ 
(resp. $\delta$). Remember that we denote $\Om_{\delta}=\Om\setminus\overline{\om}_{\delta}$. We consider the problem of finding 
$u_{\delta}\in \mH^{1}_{0}(\Omega_{\delta})$ such that 
\begin{equation}\label{NumericalPb}
-\Delta u_{\delta} = 0\quad \textrm{in}\;\;\Omega_{\delta}
\quad\textrm{and}\quad u_{\delta} = g\quad\textrm{on}\;\;\partial\Omega, 
\quad u_{\delta} = 0\quad\textrm{on}\;\;\partial\omega_{\delta}.
\end{equation}
Admittedly (\ref{NumericalPb}) is not exactly of the same form as (\ref{PbInit}). However the analysis developed in the previous 
sections can be adapted in a straightforward manner to deal with (\ref{NumericalPb}) and the results are the same. For such a 
configuration, the exact solution $u_{\delta}$ is given by
\[
\begin{array}{|ll}
u_{\delta}(\bfx) = \dsp 1 -\ln|\bfx|/\ln\delta & \mbox{ for }g=1\\[5pt]
\dsp{ u_{\delta}(\bfx) = \frac{(\delta/\vert\bfx\vert)^{-n}-(\delta/\vert\bfx\vert)^{n}}{\delta^{-n}-\delta^{n}}
\,\sin(n\theta) }
& \mbox{ for }g=\sin(n\theta),\ n\in\{1,2,\dots\}.
\end{array}
\]
Note that, with  $\omega_{\delta} = \mD_{\delta}$, we have 
$\omega = \omega_1 = \mD_{1}$ so that the logarithmic capacity potential $P$ defined by (\ref{Profiles}) verifies $P(\bfxi) = (2\pi)^{-1}\ln\vert\bfxi\vert^{-1}$. 
As a consequence, the parameter $P_{0}$ appearing in the definition of $b_{\delta}(\cdot)$ (see (\ref{decompoBis})) satisfies $P_{0} = 0$. For the computation of this parameter in other geometries, we refer the reader to \cite{Rans11}. Let us consider $\Om^h$ a polygonal approximation of the domain $\Om$. Introduce $(\mathcal{T}^h)_h$ a shape regular family of triangulations of $\overline{\Om}\!\,^h$. 
Here, $h$ denotes the average mesh size. Define the family of finite element spaces
$$
\mV^h_{\kappa} := \left\{\varphi\in \mH_{0}^{1}(\Om^h)\mbox{ such that }\varphi|_{\tau}\in \mathbb{P}_{\kappa}(\tau)\mbox{ for all }\tau\in\mathcal{T}^h\right\},
$$
where $\mathbb{P}_{\kappa}(\tau)$ is the space of polynomials of degree at most $\kappa\in\{1,2,3\}$ on the triangle $\tau$. 
We will denote $w_{\delta1}^{h}$, $w_{\delta2}^{h}$ and $w_{\delta3}^{h}$ the numerical solutions of (\ref{DiscreteVFBis}) obtained respectively with $\mV^h_{1}$, $\mV^h_{2}$
and $\mV^h_{3}$. The cut-off function $\chi$ appearing in the definition of $b_{\log}(\cdot)$ (see (\ref{DefTerVaria})) is chosen in  
$\mathscr{C}^{\infty}(\overline{\Om})$ (except for the simulation of Figure \ref{OrdersC3}) with $\chi(|\bfx|)=1$ for $|\bfx|\le 0.25$ and 
$\chi(|\bfx|)=0$ for $|\bfx|\ge 0.5$. The errors are expressed in the norms $\Vert\cdot\Vert_{\mL^{2}(\Omega\setminus\mD_{\rho})}$ 
and $\Vert\cdot\Vert_{\mH^{1}(\Omega\setminus\mD_{\rho})}$ with $\rho = 0.15$. For the computations, we use the \textit{FreeFem++}\footnote{\textit{FreeFem++}, 
\url{http://www.freefem.org/ff++/}.} software while we display the results with \textit{Matlab}\footnote{\textit{Matlab}, \url{http://www.mathworks.se/}.}.\\
\newline
\textbf{On Figures \ref{Results1}, \ref{Results2}, \ref{Results3} and \ref{Results4}}, we represent the behaviour of $\Vert u_{\delta}
-u^h_{0}\Vert_{\mL^{2}(\Omega\setminus\mD_{\rho})}$, $\Vert u_{\delta}-u^h_{0}\Vert_{\mH^{1}(\Omega\setminus\mD_{\rho})}$, $\Vert u_{\delta}
-w^h_{\delta1}-b_{\delta}(w^h_{\delta1})\bfs_{\log}\Vert_{\mL^{2}(\Omega\setminus\mD_{\rho})}$, $\Vert u_{\delta}
-w^h_{\delta1}-b_{\delta}(w^h_{\delta1})\bfs_{\log}\Vert_{\mH^{1}(\Omega\setminus\mD_{\rho})}$ with respect to the mesh size 
in logarithmic scale. Figures \ref{Results1}, \ref{Results2}, \ref{Results3} and \ref{Results4} correspond respectively to 
$\delta=10^{-1}$, $\delta=10^{-2}$, $\delta=10^{-4}$ and $\delta=10^{-10}$. Here, $u^h_{0}$ is the standard $\mrm{P}1$ approximation 
of $u_0$, the $0$ order approximation of $u_{\delta}$ defined by (\ref{LimitField}). Moreover, $w^h_{\delta1}$ is the solution 
of (\ref{DiscreteVFBis}) with $\mV^h=\mV^h_{1}$ (again $\mrm{P}1$ approximation). We take $g=1$. As predicted at the end of Section 
\ref{sectionAsyExp}, we observe that the approximation of $u_{\delta}$ by $u^h_{0}$ does not provide satisfactory results (even 
for $\delta=10^{-10}$). This is due to the error in the model, of order $|\ln\delta|^{-1}$, which decays very slowly as $\delta$ 
tends to zero. Conversely, $w^h_{\delta1}+b_{\delta}(w^h_{\delta1})\bfs_{\log}$ appears as a good approximation of $u_{\delta}$ and the rates of convergence are 
as expected. Moreover, the curves for $\delta=10^{-10}$ confirm the absence of any locking phenomenon for this numerical scheme.\\
\newline
\textbf{On Figures \ref{Orders}, \ref{OrdersC3}}, we display the behaviour of $\Vert u_{\delta}
-w^h_{\delta1}-b_{\delta}(w^h_{\delta1})\bfs_{\log}\Vert_{\mH^{1}(\Omega\setminus\mD_{\rho})}$, $\Vert u_{\delta}
-w^h_{\delta2}-b_{\delta}(w^h_{\delta2})\bfs_{\log}\Vert_{\mH^{1}(\Omega\setminus\mD_{\rho})}$, $\Vert u_{\delta}
-w^h_{\delta3}-b_{\delta}(w^h_{\delta3})\bfs_{\log}\Vert_{\mH^{1}(\Omega\setminus\mD_{\rho})}$ with respect to the 
mesh size in logarithmic scale. For the experiments of Figure \ref{Orders}, the cut-off function $\chi$ appearing in the 
definition of $b_{\log}(\cdot)$ (see (\ref{DefTerVaria})) is chosen equal to $\chi_{\mrm{exp}}$, an element of 
$\mathscr{C}^{\infty}(\overline{\Om})$ built with the exponential function. For the simulations of Figure \ref{OrdersC3}, we 
take $\chi$ equal to $\chi_{\mrm{pol}}\in\mathscr{C}^{3}(\overline{\Om})\setminus\mathscr{C}^{4}(\overline{\Om})$, a piecewise 
polynomial function of degree $7$. We take $g=1$ and $\delta=10^{-10}$. We notice that with $\chi=\chi_{\exp}$, we obtain optimal rates 
of convergence. This is not the case for $\mrm{P}3$ approximation when we choose $\chi=\chi_{\mrm{pol}}\in\mathscr{C}^{3}(
\overline{\Om})\setminus\mathscr{C}^{4}(\overline{\Om})$. However, we also remark that for the mesh sizes $h$ considered here, 
the error is smaller when $\chi$ is a polynomial function ($\chi=\chi_{\mrm{pol}}$) than when $\chi$ is built with the 
exponential function ($\chi=\chi_{\exp}$).\\
\newline
\textbf{On Figure \ref{ResultsThreshold}}, we observe the behaviour of $\Vert u_{\delta}
-w^h_{\delta2}-b_{\delta}(w^h_{\delta2})\bfs_{\log}\Vert_{\mH^{1}(\Omega\setminus\mD_{\rho})}$ with respect to the mesh 
size in logarithmic scale and for different values of $\delta$. Here, $w^h_{\delta2}$ is the solution of (\ref{DiscreteVFBis}) 
with $\mV^h=\mV^h_{2}$. We take $g=1+\sin(\theta)$. We see some thresholds in the convergence with respect to $h$: according 
to the value of $\delta$, the error stops decreasing at some $h_0$. This corresponds again to the 
error of the model. Estimates (\ref{approxFF}) indicates that this error behaves like $\delta\,\vert\ln\delta\vert$. This is 
better than the error of the $0$ order model (in $|\ln\delta|^{-1}$ ), but when $\delta$ is not so small, it is natural that 
it appears. These thresholds are absent in the curves of Figure \ref{Results1} because of the value of the source term. A 
natural approach to decrease the error consists in considering a model of higher order. Then, working as in Section 
\ref{sectionConstruction}, one can derive a model problem which does not suffer from numerical locking effect and whose solution 
yields a better approximation of $u_{\delta}$. We emphasize that at any order, this technique requires only one numerical resolution, 
and remains robust as $\delta\to 0$.

\quad\\[20pt]

\newpage

\newpage

\newpage

\section*{Acknowledgments} 
The authors would like to thank Sergey A. Nazarov, of the Faculty of Mathematics and Mechanics of St. Petersburg State University, for useful discussions and remarks. Besides, the work of the first author was supported by the Academy of Finland (decision 140998). 

\bibliography{Bibli}

\def\cprime{$'$}
\begin{thebibliography}{10}

\bibitem{MR3058835}
X.~Antoine, K.~Ramdani, and B.~Thierry.
\newblock Wide frequency band numerical approaches for multiple scattering
  problems by disks.
\newblock {\em J. Algorithms Comput. Technol.}, 6(2):241--259, 2012.

\bibitem{BaSu92}
I.~Babu{\v{s}}ka and M.~Suri.
\newblock On locking and robustness in the finite element method.
\newblock {\em SIAM J. Numer. Anal.}, 29(5):1261--1293, 1992.

\bibitem{BeBl99}
T.~Belytschko and T.~Black.
\newblock Elastic crack growth in finite elements with minimal remeshing.
\newblock {\em Int. J. Numer. Meth. Eng.}, 45(5):601--620, 1999.

\bibitem{BHBo05}
{M.F.} Ben~Hassen and E.~Bonnetier.
\newblock Asymptotic formulas for the voltage potential in a composite medium
  containing close or touching disks of small diameter.
\newblock {\em {Multiscale Model. Simul.}}, 4(1):250--277, 2005.

\bibitem{865332}
J.-P. B\'{e}renger.
\newblock A multiwire formalism for the {FDTD} method.
\newblock {\em IEEE Trans. Electromagn. Compat.on}, 42(3):257--264, 2000.

\bibitem{BoDa13}
V.~Bonnaillie-No{\"e}l and M.~Dambrine.
\newblock Interactions between moderately close circular inclusions: the
  {D}irichlet-{L}aplace equation in the plane.
\newblock {\em Asymptot. Anal.}, 84(3-4):197--227, 2013.

\bibitem{BDTV07}
V.~Bonnaillie-No{\"e}l, M.~Dambrine, S.~Tordeux, and G.~Vial.
\newblock {On moderately close inclusions for the Laplace equation}.
\newblock {\em C. R. Acad. Sci., Ser. I}, 345(11):609--614, 2007.

\bibitem{BDTV09}
V.~{Bonnaillie-No{\"e}l}, M.~Dambrine, S.~Tordeux, and G.~Vial.
\newblock Interactions between moderately close inclusions for the {L}aplace
  equation.
\newblock {\em Math. Models Methods Appl. Sci.}, 19(10):1853--1882, 2009.

\bibitem{BoVo00}
E.~{Bonnetier} and M.~{Vogelius}.
\newblock {An elliptic regularity result for a composite medium with
  ''touching'' fibers of circular cross-section.}
\newblock {\em {SIAM J. Math. Anal.}}, 31(3):651--677, 2000.

\bibitem{BDLN92}
M.~{Bourlard}, M.~{Dauge}, {M.-S.} {Lubuma}, and S.~{Nicaise}.
\newblock {Coefficients of the singularities for elliptic boundary value
  problems on domains with conical points. III: Finite element methods on
  polygonal domains}.
\newblock {\em {SIAM J. Numer. Anal.}}, 29(1):136--155, 1992.

\bibitem{BrSc08}
S.C. Brenner and L.R Scott.
\newblock {\em {The mathematical theory of finite element methods. 3rd ed.}}
\newblock Springer, New York, 2008.

\bibitem{MR1429996}
A.~Campbell and {S.A.} Nazarov.
\newblock Une justification de la m\'ethode de raccordement des
  d\'eveloppements asymptotiques appliqu\'ee \`a un probl\`eme de plaque en
  flexion. {E}stimation de la matrice d'imp\'edance.
\newblock {\em J. Math. Pures Appl.}, 76(1):15--54, 1997.

\bibitem{MR2991243}
M.~Cassier and C.~Hazard.
\newblock Multiple scattering of acoustic waves by small sound-soft obstacles
  in two dimensions: mathematical justification of the {F}oldy-{L}ax model.
\newblock {\em Wave Motion}, 50(1):18--28, 2013.

\bibitem{MR1990198}
Z.~Chen and X.~Yue.
\newblock Numerical homogenization of well singularities in the flow transport
  through heterogeneous porous media.
\newblock {\em Multiscale Model. Simul.}, 1(2):260--303, 2003.

\bibitem{CJKLZ05a}
P.~{Ciarlet Jr.}, B.~Jung, S.~Kaddouri, S.~Labrunie, and J.~Zou.
\newblock {The Fourier singular complement method for the Poisson problem. Part
  I: Prismatic domains}.
\newblock {\em Numer. Math.}, 101(3):423--450, 2005.

\bibitem{CJKLZ05b}
P.~{Ciarlet Jr.}, B.~Jung, S.~Kaddouri, S.~Labrunie, and J.~Zou.
\newblock {The Fourier Singular Complement Method for the Poisson problem. Part
  II: axisymmetric domains}.
\newblock {\em Numer. Math.}, 102(4):583--610, 2006.

\bibitem{CiLa11}
P.~{Ciarlet Jr.} and S.~Labrunie.
\newblock Numerical solution of maxwell's equations in axisymmetric domains
  with the fourier singular complement method.
\newblock {\em J. Differ. Equ. Appl.}, 3:113--155, 2011.

\bibitem{DaVi05}
M.~{Dambrine} and G.~{Vial}.
\newblock {Influence of a boundary perforation on the Dirichlet energy}.
\newblock {\em {Control Cybern.}}, 34(1):117--136, 2005.

\bibitem{DaVi07}
M.~Dambrine and G.~Vial.
\newblock A multiscale correction method for local singular perturbations of
  the boundary.
\newblock {\em Math. Mod. Num. Anal.}, 41(01):111--127, 2007.

\bibitem{DoBe99}
J.~Dolbow and T.~Belytschko.
\newblock A finite element method for crack growth without remeshing.
\newblock {\em Int. J. Numer. Meth. Eng.}, 46(1):131--150, 1999.

\bibitem{Duar96}
{C.A.} Duarte and {J.T.} Oden.
\newblock An \textit{h}-\textit{p} adaptive method using clouds.
\newblock {\em Comput. Methods in Appl. Mech. Eng.}, 139(1):237--262, 1996.

\bibitem{MR2477579}
Y.~Efendiev and T.Y. Hou.
\newblock {\em Multiscale finite element methods}, volume~4 of {\em Surveys and
  Tutorials in the Applied Mathematical Sciences}.
\newblock Springer, New York, 2009.
\newblock Theory and applications.

\bibitem{HaLo02}
C.~Hazard and S.~Lohrengel.
\newblock A singular field method for maxwell's equations: Numerical aspects
  for 2d magnetostatics.
\newblock {\em {SIAM J. Numer. Anal.}}, 40(3):1021--1040, 2002.

\bibitem{MR630701}
{A.M.} {Il'in}.
\newblock Study of the asymptotic behavior of the solution of an elliptic
  boundary value problem in a domain with a small hole.
\newblock {\em Trudy Sem. Petrovsk.}, (6):57--82, 1981.

\bibitem{Ilin92}
{A.M.} {Il'in}.
\newblock {\em Matching of asymptotic expansions of solutions of boundary value
  problems}, volume 102 of {\em Translations of Mathematical Monographs}.
\newblock AMS, Providence, RI, 1992.

\bibitem{Kond67}
{V.A.} {Kondratiev}.
\newblock Boundary-value problems for elliptic equations in domains with
  conical or angular points.
\newblock {\em Trans. Moscow Math. Soc.}, 16:227--313, 1967.

\bibitem{KoMR97}
{V.A.} {Kozlov}, V.G. {Maz'ya}, and J.~{Rossmann}.
\newblock {\em Elliptic Boundary Value Problems in Domains with Point
  Singularities}, volume~52 of {\em Mathematical Surveys and Monographs}.
\newblock AMS, Providence, 1997.

\bibitem{MR0350075}
N.~N. Lebedev.
\newblock {\em Special functions and their applications}.
\newblock Dover Publications, Inc., New York, 1972.

\bibitem{MR912054}
{V.G.} {Maz'ya} and {S.A.} Nazarov.
\newblock Asymptotic behavior of energy integrals under small perturbations of
  the boundary near corner and conic points.
\newblock {\em Trudy Moskov. Mat. Obshch.}, 50:79--129, 1987.

\bibitem{Na61}
{V.G.} {Maz'ya}, {S.A.} Nazarov, and {B.A.} Plamenevski{\u\i}.
\newblock Asymptotic expansions of eigenvalues of boundary value problems for
  the {L}aplace operator in domains with small openings.
\newblock {\em Izv. Akad. Nauk SSSR Ser. Mat.}, 48(2):347--371, 1984.

\bibitem{MaNP00}
{V.G.} {Maz'ya}, {S.A.} Nazarov, and {B.A.} Plamenevski{\u\i}.
\newblock {\em {Asymptotic theory of elliptic boundary value problems in
  singularly perturbed domains, Vol. 1, 2}}.
\newblock {Birkh\"{a}user}, Basel, 2000.
\newblock Translated from the original German 1991 edition.

\bibitem{MeBa96}
{J.M.} Melenk and I.~Babu{\v{s}}ka.
\newblock The partition of unity finite element method: basic theory and
  applications.
\newblock {\em Comput. Methods in Appl. Mech. Eng.}, 139(1-4):289--314, 1996.

\bibitem{Naza99}
{S.A.} Nazarov.
\newblock Asymptotic conditions at a point, selfadjoint extensions of
  operators, and the method of matched asymptotic expansions.
\newblock In {\em Proceedings of the {S}t. {P}etersburg {M}athematical
  {S}ociety, {V}ol. {V}}, volume 193 of {\em Amer. Math. Soc., Transl. Ser. 2},
  pages 77--125, Providence, RI, 1999.

\bibitem{NaSo03}
{S.A.} Nazarov and J.~Soko{\l}owski.
\newblock Asymptotic analysis of shape functionals.
\newblock {\em J. Math. Pures Appl.}, 82(2):125--196, 2003.

\bibitem{NaSo06}
{S.A.} Nazarov and J.~Soko{\l}owski.
\newblock Self-adjoint extensions for the {N}eumann {L}aplacian and
  applications.
\newblock {\em Acta Math. Sin. (Engl. Ser.)}, 22(3):879--906, 2006.

\bibitem{Peaceman}
D.W. Peaceman.
\newblock Interpretation of well-block pressures in numerical reservoir
  simulations.
\newblock {\em Soc. Pet. Eng. J.}, 18(3):183--194, 1978.

\bibitem{1353491}
C.J. Railton, B.P. Koh, and I.J. Craddock.
\newblock The treatment of thin wires in the {FDTD} method using a weighted
  residuals approach.
\newblock {\em IEEE Trans. Antennas Propag.}, 52(11):2941--2949, 2004.

\bibitem{Rans11}
T.~Ransford.
\newblock Computation of logarithmic capacity.
\newblock {\em Comput. Meth. Funct. Theor.}, 10(2):555--578, 2011.

\bibitem{Scot76}
R.~Scott.
\newblock Optimal {$L^{\infty}$} estimates for the finite element method on
  irregular meshes.
\newblock {\em Math. Comp.}, 30(136):681--697, 1976.

\end{thebibliography}
\bibliographystyle{plain}
\end{document}